\documentclass[aos,nameyear,seceqn,dvips]{arximspdf}
\usepackage{graphicx}

\doi{10.1214/10-AOS846}
\volume{39}
\issue{1}
\pubyear{2011}
\firstpage{389}
\lastpage{417}

\makeatletter

\renewcommand{\epsilon}{\varepsilon}

\newcommand{\vb}{{\mathbf b}}
\newcommand{\vc}{{\mathbf c}}
\newcommand{\ve}{{\mathbf e}}
\newcommand{\vS}{{\mathbf S}}
\newcommand{\vE}{{\mathbf E}}
\newcommand{\vG}{{\mathbf G}}
\newcommand{\vF}{{\mathbf F}}
\newcommand{\vu}{{\mathbf u}}
\newcommand{\vA}{{\mathbf A}}
\newcommand{\vB}{{\mathbf B}}
\newcommand{\vC}{{\mathbf C}}
\newcommand{\vH}{{\mathbf H}}
\newcommand{\vI}{{\mathbf I}}
\newcommand{\vL}{{\mathbf L}}
\newcommand{\vM}{{\mathbf M}}
\newcommand{\vN}{{\mathbf N}}
\newcommand{\vQ}{{\mathbf Q}}
\newcommand{\vR}{{\mathbf R}}
\newcommand{\vV}{{\mathbf V}}
\newcommand{\vY}{{\mathbf Y}}
\newcommand{\vX}{{\mathbf X}}
\newcommand{\vD}{{\mathbf D}}

\newcommand{\veins}{{\mathbf 1}}
\newcommand{\vnull}{{\mathbf 0}}

\newcommand{\vSigma}{\bolds\Sigma}
\newcommand{\valpha}{\bolds\alpha}
\newcommand{\vbeta}{\bolds\beta}
\newcommand{\vtau}{\bolds\tau}
\newcommand{\vmu}{\bolds\mu}
\newcommand{\vpi}{\bolds\pi}

\newcommand{\E}{\mathrm{E}}
\newcommand{\Var}{\operatorname{Var}}
\newcommand{\Cov}{\operatorname{Cov}}
\newcommand{\diag}{\operatorname{diag}}

\newcommand{\ra}{\rightarrow}

\newtheorem{lemma}{Lemma}[section]
\newtheorem{theorem}[lemma]{Theorem}
\newtheorem{corollary}[lemma]{Corollary}
\newproclaim{example}{Example}
\newproclaim{remark}{Remark}
\newcommand{\vepsilon}{\bolds\epsilon}

\makeatother

\begin{document}
\begin{frontmatter}

\title{GEE analysis of clustered binary data with diverging number of covariates}
\runtitle{GEE with large $p$}

\begin{aug}
\author[A]{\fnms{Lan} \snm{Wang}\thanksref{t1}\corref{}\ead[label=e1]{lan@stat.umn.edu}}
\thankstext{t1}{Supported by NSF Grant DMS-1007603. }
\runauthor{L. Wang}
\affiliation{University of Minnesota}
\address[A]{School of Statistics\\
University of Minnesota\\
224 Church Street, SE\\
Minneapolis, Minnesota 55455\\
USA\\
\printead{e1}} 
\end{aug}

\received{\smonth{12} \syear{2009}}
\revised{\smonth{7} \syear{2010}}

\begin{abstract}
Clustered binary data with a large number of covariates have become
increasingly  common in many scientific disciplines. This paper develops
an asymptotic theory for generalized estimating equations (GEE) analysis of
clustered binary data when the number of covariates grows to infinity with the
number of clusters. In this ``large $n$, diverging $p$'' framework, we provide
appropriate regularity conditions and establish the existence, consistency and
asymptotic normality of the GEE estimator. Furthermore, we prove that the
sandwich variance formula remains valid. Even when the working correlation
matrix is misspecified, the use of the sandwich variance formula leads to
an asymptotically valid confidence interval and Wald test for an estimable linear
combination of the unknown parameters. The accuracy of the asymptotic
approximation is examined via numerical simulations. We also discuss the
``diverging $p$'' asymptotic theory for general GEE. The results in this paper
extend the recent elegant work of Xie and Yang [\textit{Ann. Statist.} \textbf{31} (2003) 310--347] and Balan and
Schiopu-Kratina [\textit{Ann. Statist.} \textbf{32} (2005) 522--541] in the ``fixed $p$'' setting.
\end{abstract}



\end{frontmatter}

\section{Introduction}
A fundamental problem in statistical analysis is to characterize the effects
of a set of covariates $X_1,\ldots,X_p$ on a response variable $Y$ based on a
sample of size $n$. Recently, there has been considerable interest in
investigating this problem in the so-called ``large $n$, diverging
$p$''
asymptotic framework, where the dimension of the covariates increases to
infinity with the sample size. This setup allows statisticians to adopt a more
complex statistical model as more abundant data become available, and thus to
reduce the modeling bias.

The ``large $n$, diverging $p$'' framework can be traced back to the earlier
pioneering work on M-estimators with a diverging number of parameter; see Huber
(\citeyear{hub1973}), Portnoy (\citeyear{por1984},
\citeyear{por1985}, \citeyear{por1988}),
Mammen (\citeyear{mam1989}),
Welsh (\citeyear{wel1989}), Bai and Wu (\citeyear{baiwu1994}), He and Shao (\citeyear{hesha2000}) and the references
therein. With the advent of high-dimensional data in many scientific areas,
statistical theory developed in this new framework has become crucial for
guiding practical data analysis with high-dimensional covariates, which relies
heavily on asymptotic theory to justify its validity. By allowing the
covariates' dimension to increase with the sample size, Fan and Peng (\citeyear{fanpen2004})
studied nonconcave penalized likelihood; Lam and Fan (\citeyear{lamfan2008}) investigated
profile-kernel likelihood inference with generalized varying coefficient
partially linear models; Huang, Horowitz and Ma (\citeyear{huahorma2008}) explored bridge
estimators in linear regression; Hjort, McKeague  and Van Keilegom (\citeyear{hjomckvan2008}) and
Chen, Peng and Qin (\citeyear{chepenqin2009}) studied the effects of data dimension on empirical
likelihood; Zou and Zhang (\citeyear{zouzha2009}) studied the adaptive elastic net, Zhu and Zhu~(\citeyear{zhuzhu2009})
investigated parameter estimation in a semiparametric regression model
with highly correlated predictors. In the aforementioned
literature, the number of covariates $p$ grows to infinity at a polynomial
rate $o(n^{\alpha})$ for some $0<\alpha<1$. In particular, most of these
papers provide necessary conditions under which classical asymptotic
theories remain valid for $\alpha$ in the range $[\frac{1}{5},\frac{1}{2}]$.

A different line of research considers the case where $p$ can be much larger than
$n$ and even grow at an exponential rate of $n$, in which case the sparsity
assumption and other more stringent regularity conditions are generally
required to investigate the large-sample properties. Furthermore, it is worth
noting that much work has also been devoted to classification and
multiple hypotheses testing problems with high-dimensional covariates, but
these problems are different in nature from what is discussed in this paper.
We refer to the review papers of Donoho (\citeyear{don2000}), Fan and Li (\citeyear{fanli2006}) and Fan and
Lv (\citeyear{fanlv2009}) for more comprehensive references on high-dimensional data analysis.

When the research focus is on modeling the relationship between $Y$ and a
high-dimensional vector of covariates, the existing literature in the ``large
$n$, diverging $p$'' setting has been largely restricted to independent data.
In many modern data sets, in addition to the large dimensionality of
covariates, complexity also arises when the responses are correlated due to
repeated measures or clustered design. One representative example is the
Framingham Heart Study, where the researchers are interested in linking common
risk factors to the occurrence of cardiovascular diseases. In this study,
many variables, such as age, smoking status, cholesterol level and blood
pressure, were recorded for the participants during their clinic visits
over the years to describe their physical characteristics and lifestyles.
Another example is the Chicago Longitudinal Study in social science, which
investigated the educational and social development of about 1500 low income,
minority youths in the Chicago area. The study collected a large amount of
information on many variables that measure children's early antisocial
behavior, individual-level attributes of the child, family attributes and
social characteristics of both the child and the family, among others. In some
other examples of clustered data, the number of variables measured for each
individual or experimental unit may not be many, but when one considers
various interaction effects, the actual number of predictors in the statistical model can
still be large and better fits the ``large $p$'' setup.

The intrinsic complexity of clustered data raises challenging issues for
statistical analysis, especially for correlated non-Gaussian data where it is
difficult to specify the full likelihood. In this paper, we establish the
asymptotic properties of generalized estimating equations (GEE), a
semiparametric procedure widely used in practice for clustered data analysis,
while allowing the covariate dimension to grow to infinity with the sample
size.

The GEE procedure was introduced in a seminal paper of Liang and\break  Zeger~(\citeyear{liazeg1986})
as a useful extension of generalized linear models [McCullagh and Nelder
(\citeyear{mccnel1989})]
to correlated data. Instead of specifying the full likelihood, it only
requires the knowledge of the first two marginal moments and a working
correlation matrix. Thus, it is particularly effective for modeling clustered
binary or count data. A key advantage of the GEE approach is that it yields a
consistent estimator (in the classical ``large $n$, fixed $p$'' setup), even if
the working correlation structure is misspecified. The GEE estimator is also
asymptotically efficient if the correlation structure is indeed correctly
specified. The original paper of Liang and Zeger focused mostly on the
methodology development. Li (\citeyear{li1997}) adopted a minimax approach to study the
consistency of GEE. A more complete and systematic large-sample theory
for GEE, including consistency and asymptotic normality, was elegantly
established by Xie and Yang (\citeyear{xieyan2003}). Balan and Schiopu-Kratina (\citeyear{balsch2005}) also
rigorously studied a closely related pseudo-likelihood framework for GEE.
However, these papers all assume that $p$ is fixed and that the number of clusters
$n$ goes to infinity. Xie and Yang (\citeyear{xieyan2003}) also considered the case where the cluster
size (number of observations within each cluster) is itself  large, which
corresponds to a large number of time points in the longitudinal
setting.\looseness=-1

This paper examines the effect of high-dimensional covariates on the GEE
estimator in the ``large $n$, diverging $p$'' setup, where $p=p_n$ is a
function of the sample size $n$. We focus on clustered binary data because
binary response (e.g., disease status) is ubiquitous in many scientific
applications and because of the relative transparency of technical derivation. We also
discuss the related theory for general GEE in Section \ref{sec51} The main technical
challenges come from the high dimensionality of the covariates, the
dependence among observations within each cluster and the nuisance parameters
in the working correlation matrix. We provide a self-contained derivation and
extend earlier theory in the literature on M-estimation with a large number of
parameters, which is not tailored for clustered data and generally has not
considered nuisance parameters.

We aim to answer the following essential questions. To what extent can the asymptotic
results derived in the classical asymptotic framework for GEE still be deemed
trustworthy when the number of covariates is large? How large can $p_n$ be
(relative to $n$)? The main findings in this paper reveal that under
reasonable conditions, the GEE estimator\vspace*{1pt} $\widehat{\vbeta}_n$ is
$\sqrt{p_n/n}$-consistent when $p_n^2/n\ra 0$ and that an arbitrary linear
combination $\valpha_n^{T}(\widehat{\vbeta}_n-\vbeta_{n0})$ is asymptotically
normal when $p_n^3/n\ra 0$, where $\vbeta_{n0}$ is the true parameter value.
These findings resonate with those in the literature for independent data in
the ``large $p$'' setting. Moreover, we also verify that the desirable
robustness property against working correlation matrix misspecification still
holds and that both the sandwich variance formula and the large-sample Wald
test still remain valid in this new context. Understanding these fundamental
questions is essential to justifying asymptotic statistical inference based on
GEE  for analyzing real-world clustered data containing many covariates, such
as the validity of the confidence intervals provided by the GEE package in R,
SAS and other statistical software packages.

The rest of the paper is organized as follows. In Section \ref{sec2}, we provide a
brief review of the GEE procedure for analyzing clustered binary data. Section
\ref{sec3} establishes the consistency and asymptotic normality of the GEE estimator,
the consistency of the sandwich variance formula and the validity of the
large-sample Wald test in the ``large $n$, diverging $p$'' framework. Section
\ref{sec4}
examines the asymptotic results via numerical simulations. Section \ref{sec5} discusses
general GEE and related problems.

\section{Generalized estimating equations}\label{sec2}
For the $j$th observation of the $i$th cluster, we observe a binary response
variable $Y_{ij}$ and a $p_n$-dimensional vector of covariates $\vX_{ij}$,
$i=1,\ldots,n$ and $j=1,\ldots,m_i$. Observations from different clusters are
independent, but those from the same clusters are correlated. Let
$\vY_i=(Y_{i1},\ldots,Y_{im_i})^T$ denote the vector of responses for the
$i$th cluster and let $\vX_i=(\vX_{i1}, \ldots,\vX_{im_i})^T$ be the
associated $m_i\times p_n$ matrix of covariates.

The marginal regression approach of GEE assumes that
$\E(Y_{ij}|\vX_{ij})=\pi_{ij}$ and
$\Var(Y_{ij}|\vX_{ij})=\pi_{ij}(1-\pi_{ij})$, where a dispersion parameter may
be added in the marginal variance function if overdispersion is suspected to
be present. Furthermore, it relates the covariates to the marginal mean by
specifying that
\begin{equation}\label{marginal}
\operatorname{logit}(\pi_{ij})=\vX_{ij}^T\vbeta_n,
\end{equation}
where
$\operatorname{logit}(\pi_{ij})=\log(\frac{\pi_{ij}}{1-\pi_{ij}})$ is the link
function and $\vbeta_n$ is a $p_n$-dimensional vector of parameters. The true
unknown parameter value is denoted by $\vbeta_{n0}$.

Let $\vpi_i(\vbeta_n)=(\pi_{i1}(\vbeta_n),\ldots,\pi_{im_i}(\vbeta_n))^T$,
where
$\pi_{ij}(\vbeta_n)=\exp(\vX_{ij}^T\vbeta_n)/[1+\exp(\vX_{ij}^T\vbeta_n)]$.
Further, let $\vA_i(\vbeta_n)$ be the $m_i \times m_i$ diagonal matrix with
the $j$th diagonal element
$\vA_{ij}(\vbeta_n)=\pi_{ij}(\vbeta_n)(1-\pi_{ij}(\vbeta_n))$,
$j=1,\ldots,m_i$.
In what follows, we assume $m_i= m<\infty,$ for simplicity. Liang and Zeger
(\citeyear{liazeg1986}) suggested to estimate $\vbeta_{n0}$ by solving the following
generalized estimating equation in $\vbeta_n$:
\begin{equation}\label{GEE1}
\sum_{i=1}^n\vX_i^T\vA_i(\vbeta_n)\vV_i^{-1}\bigl(\vY_i-\vpi_i(\vbeta_n)\bigr)=0,
\end{equation}
where $\vV_i$ is a working covariance matrix.

\section{Asymptotic properties when $p_n\ra \infty$}\label{sec3}
\subsection{GEE estimator with estimated working correlation
matrix}\label{sec31}

In applications, the true correlation matrix of $\vY_i$, denoted by $\vR_0$, is
unknown. The working covariance matrix is often specified via a working
correlation matrix $\vR(\vtau)$:
$\vV_i=\vA_i^{1/2}(\vbeta_n)\vR(\vtau)\vA_i^{1/2}(\vbeta_n),$ where $\vtau$ is
a finite-dimensional\vadjust{\goodbreak} parameter. Commonly used working correlation
structures include $\mbox{AR-1}$, compound symmetry and unstructured working
correlation, among others. Note that, in practice, the working correlation matrix is chosen to be independent
 of the covariates, for simplicity. However, for correlated non-normal data, the range of correlation generally depends
on the univariate marginals. Thus, $ \vR(\vtau)$ should be understood as a weight matrix
[Chaganty and Joe (\citeyear{chajoe2004})].
Chaganty and Joe demonstrated that GEE with an appropriately chosen working
correlation matrix does have good efficiency when compared with a proper likelihood model.

Given a working correlation structure, $\vtau$ is
often estimated using a residual-based moment method, which requires an initial
consistent estimator of $\vbeta_{n0}$. We use $\widehat{\vR}$ to denote the
resulting estimated working correlation matrix, with the subscript
``$n$''
suppressed. Following (\ref{GEE1}), we formally define the \textit{GEE estimator}
$\widehat{\vbeta}_n$ as the solution of
\begin{equation}
\vS_n(\vbeta_n)=\sum_{i=1}^n\vX_i^T\vA_i^{1/2}
(\vbeta_n)\widehat{\vR}^{-1}\vA_i^{-1/2}(\vbeta_n)\bigl(\vY_i-\vpi_i(\vbeta_n)\bigr)=0.
\end{equation}
To solve for $\widehat{\vbeta}_n$, we can iterate between a modified
Fisher scoring algorithm for $\vbeta_n$ and the moment estimation for $\vtau$.
In the following, we provide examples of an initial consistent estimator and
an estimated working correlation matrix.

\begin{example}[(Initial estimator for $\vbeta_{n0}$ when $p_n\ra
\infty$)]\label{exam1}
A simple way to obtain an initial estimator for $\vbeta_{n0}$ is
to solve the generalized estimating equations under the working independence
assumption
\begin{eqnarray}\label{IEE}
\widetilde{S}_n(\vbeta_n)=\sum_{i=1}^n\vX_i^T\bigl(\vY_i-\vpi_i(\vbeta_n)\bigr)=0.
\end{eqnarray}

Under conditions (A1)--(A3) in Section \ref{sec32}, we can show that if $p_n^2/n\ra 0$
as $n\ra \infty$, then the independence estimating equations in (\ref{IEE})
have a root $\widetilde{\vbeta}_n$ such that
\begin{eqnarray}\label{inia}
\Vert\widetilde{\vbeta}_n-\vbeta_{n0}\Vert=O_p\bigl(\sqrt{p_n/n}\bigr),
\end{eqnarray}
where
$\Vert\cdot\Vert$ denotes the Euclidean norm of a vector. A detailed derivation of
(\ref{inia}) is given in the \hyperref[app]{Appendix}.
\end{example}

\begin{example}[(Estimated working correlation matrix when $p_n\ra
\infty$)]\label{exam2}
In Balan and Schiopu-Kratina (\citeyear{balsch2005}), it was suggested to use
\[
\widehat{\vR}=\frac{1}{n}\sum_{i=1}^n\vA_i^{-1/2}(\widetilde{\vbeta}_n)\bigl(\vY_i-\vpi_i(\widetilde{\vbeta}_n)\bigr)
\bigl(\vY_i-\vpi_i(\widetilde{\vbeta}_n)\bigr)^T\vA_i^{-1/2}(\widetilde{\vbeta}_n),
\]
where $\widetilde{\vbeta}_n$ is a preliminary $\sqrt{n/p_n}$-consistent estimator of $\vbeta_{n0}$, such
as the one discussed in Example \ref{exam1}. This provides a moment estimator of the
unstructured working correlation matrix. Assuming conditions (A1)--(A3) of
Section \ref{sec32}, we can prove that if $p_n^2/n\ra 0$ as $n\ra \infty$, then
%
\begin{equation}\label{inib}
\Vert\widehat{\vR}^{-1}-\vR_0^{-1}\Vert=O_p\bigl(\sqrt{p_n/n}\bigr),
\end{equation}
where $\vR_0$ denotes the true common correlation matrix. Here, and throughout the paper,
for a matrix $\vB$, $\Vert\vB\Vert=[\operatorname{Tr}(\vB\vB^T)]^{1/2}$ denotes its
Frobenius norm. A detailed derivation of (\ref{inib}) is given in the supplementary article [Wang (\citeyear{wan2010})].
\end{example}

\subsection{Existence and consistency}\label{sec32}
In Fan and Peng (\citeyear{fanpen2004}), Lam and Fan (\citeyear{lamfan2008}) and Huang, Horowitz and Ma (\citeyear{huahorma2008}), the
estimator is defined as the minimizer of a certain objective function. We use
alternative techniques here to establish the existence and consistency of the
GEE estimator, which involve the roots of estimating equations. The
approach we adopt here is also different from that of Xie and Yang~(\citeyear{xieyan2003}) and
Balan and Schiopu-Kratina (\citeyear{balsch2005}), both of which rely on properties of injective
functions.

We directly verify the following condition: $\forall\epsilon>0$, there
exists a constant $\Delta>0$ such that for all $n$ sufficiently large,
%
\begin{equation}\label{bull}
P\Bigl(\sup_{\Vert\vbeta_n-\vbeta_{n0}\Vert=
\Delta\sqrt{p_n/n}}(\vbeta_n-\vbeta_{n0})^T\vS_n(\vbeta_n)<0\Bigr)\geq
1-\epsilon.
\end{equation}
Condition (\ref{bull}) is sufficient to ensure the existence
of a sequence of roots $\widehat{\vbeta}_n$ of the equation
$\vS_n(\vbeta_n)=0$ such that
$\Vert\widehat{\vbeta}_n-\vbeta_{n0}\Vert=O_P(\sqrt{p_n/n})$. This approach follows from Theorem 6.3.4 of
Ortega and Rheinboldt (\citeyear{ortrhe1970}).  In Portnoy (\citeyear{por1984}), this technique was applied
to establish the existence and consistency of an M-estimator for i.i.d. data;
in a different setting, it was used by Wang et al. (\citeyear{wanxuezhucho2010}) to study a partial
linear single-index model. This leads to a more
straightforward and elegant proof of weak consistency. On the other hand, the method relying on injective
functions [Xie and Yang (\citeyear{xieyan2003});
Balan and Schiopu-Kratina (\citeyear{balsch2005})] can also be used to prove strong consistency.

To prove consistency and asymptotic normality, we need the following general
regularity conditions:
\begin{longlist}
\item[(A1)] $\sup_{i,j}\Vert\vX_{ij}\Vert=O(\sqrt{p_n})$;
\item[(A2)] the unknown parameter $\vbeta_n$ belongs to a compact subset
$\mathcal{B}\subseteq R^{p_n}$, the true parameter value $\vbeta_{n0}$ lies
in the interior of $\mathcal{B}$ and there exist two positive constants, $b_1$ and
$b_2$, such that $0<b_1\leq \pi_{ij}(\vbeta_{n0})\leq b_2<1$, $\forall i, j$;
\item[(A3)] there exist two positive constants, $b_3$ and $b_4,$ such that
\[
b_3\leq
\lambda_{\min}\Biggl(n^{-1}\sum_{i=1}^n\vX_i^T\vX_i\Biggr)\leq
\lambda_{\max}\Biggl(n^{-1}\sum_{i=1}^n\vX_i^T\vX_i\Biggr) \leq b_4,
\]
where
$\lambda_{\min}$ (resp. $\lambda_{\max}$) denotes the minimum (resp. maximum)
eigenvalue of a matrix;
\item[(A4)] the common true correlation matrix $\vR_0$ has eigenvalues bounded away from
zero and $+\infty$; the estimated working correlation matrix $\widehat{\vR}$
satisfies $\Vert\widehat{\vR}^{-1}-\overline{\vR}^{-1}\Vert=O_p(\sqrt{p_n/n})$,
where $\overline{\vR}$ is a constant positive definite matrix with eigenvalues
bounded away from zero and $+\infty$; we do not require $\overline{\vR}$ to be
the true correlation matrix $\vR_0$.
\end{longlist}

\begin{remark}
Condition (A1) is a common assumption in the literature
on M-estimators with diverging dimension. For example, it is the same as
assumption (3.9) of Portnoy (\citeyear{por1985}) and it is implied by conditions (C.9) and
(C.10) of Welsh~(\citeyear{wel1989}). This condition holds almost surely under some weak
moment conditions for $X_{ij}$ from spherically symmetric distributions [see,
e.g., the discussions in He and Shao (\citeyear{hesha2000})]. When $m=1$ (i.e., each
cluster has only one observation), condition (A3) is also popularly adopted in
the literature on high-dimensional regression for independent data. It can be
shown that condition (A3) is implied by the following slightly stronger
condition: there exist two positive constants, $c_1\leq c_2$, such that $\forall
1\leq j\leq m$,
\[
c_1\leq
\lambda_{\min}\Biggl(n^{-1}\sum_{i=1}^n\vX_{ij}\vX_{ij}^T\Biggr)\leq
\lambda_{\max}\Biggl(n^{-1}\sum_{i=1}^n\vX_{ij}\vX_{ij}^T\Biggr) \leq c_2.
\]
Finally, condition (A4) is a direct extension of a similar assumption in the
``fixed $p$'' case. Liang and Zeger (\citeyear{liazeg1986}) assumes that the estimator of the
working correlation matrix parameter $\widehat{\vtau}$ satisfies
$\sqrt{n}(\widehat{\vtau}-\vtau_0)=O_p(1)$ for some $\tau_0$.  Assumption (C2)
of Chen and Jin (\citeyear{chejin2006}) is of similar nature, while Xie and Yang (\citeyear{xieyan2003}) assumes
the nuisance parameter $\vtau$ to be completely known. Note
that Example \ref{exam2} in Section~\ref{sec31} guarantees that (A4) is satisfied when a
nonparametric moment estimator is used for the working correlation matrix, in
which case
$\overline{\vR}=\vR_0$.
\end{remark}

We use notation similar to that in Xie and Yang (\citeyear{xieyan2003}) and Balan and
Schiopu-Kratina (\citeyear{balsch2005}). Consider the following estimating equation:
\[
\overline{\vS}_n(\vbeta_n)=
\sum_{i=1}^n\vX_i^T\vA_i^{1/2}(\vbeta_n)\overline{\vR}^{-1}\vA_i^{-1/2}(\vbeta_n)\bigl(\vY_i-\vpi_i(\vbeta_n)\bigr).
\]
If we let $\overline{\vM}_n(\vbeta_n)$ denote the covariance matrix of
$\overline{\vS}_n(\vbeta_n)$, then
\[
\overline{\vM}_n(\vbeta_n)
=\sum_{i=1}^n\vX_i^T\vA_i^{1/2}(\vbeta_n)\overline{\vR}^{-1}\vR_0\overline{\vR}^{-1}\vA_i^{1/2}(\vbeta_n)\vX_i.
\]

To prove the consistency, the essential idea is to approximate
$\vS_n(\vbeta_n)$ by $\overline{\vS}_n(\vbeta_n)$, whose moments are easier to
evaluate. Lemma \ref{tiger1} below establishes the accuracy of this
approximation, which also plays an important role in deriving the asymptotic
normality in Section \ref{sec33}.

\begin{lemma}\label{tiger1}
Assume conditions \textup{(A1)}--\textup{(A4)}. If $n^{-1}p_n^2=o(1)$, then
\[
\Vert\vS_n(\vbeta_{n0})-\overline{\vS}_n(\vbeta_{n0})\Vert=O_p(p_n).
\]
\end{lemma}

To facilitate the Taylor expansion of the estimating function
$\vS_n(\vbeta_n)$, we also use
$\overline{\vD}_n(\vbeta_n)=-\frac{\partial}{\partial
\vbeta_{n}^T}\overline{\vS}_n(\vbeta_n)$ to approximate the negative gradient
function $\vD_n(\vbeta_n)=-\frac{\partial}{\partial
\vbeta_{n}^T}\vS_n(\vbeta_n)$. Lemma \ref{tiger2} below provides a useful
representation of $\overline{\vD}_n(\vbeta_n)$, based on which, Lemma
\ref{tiger3}
establishes the approximation of gradient functions.

\begin{lemma} \label{tiger2}
\begin{equation}\label{rep}
\overline{\vD}_n(\vbeta_n)=\overline{\vH}_n(\vbeta_n)+\overline{\vE}_n(\vbeta_n)+\overline{\vG}_n(\vbeta_n),
\end{equation}
where
\begin{eqnarray*} \overline{\vH}_n(\vbeta_n)&=&
\sum_{i=1}^n\vX_i^T\vA_i^{1/2}(\vbeta_n)\overline{\vR}^{-1}\vA_i^{1/2}(\vbeta_n)\vX_i,\\
\overline{\vE}_n(\vbeta_n)&=&\frac{1}{2}\sum_{i=1}^n\sum_{j=1}^m\bigl(1-2\pi_{ij}(\vbeta_n)\bigr)\epsilon_{ij}(\vbeta_n)
\vX_i^T\vA_i^{1/2}(\vbeta_n)\overline{\vR}^{-1}\ve_j\ve_j^T\vX_i,\\
\overline{\vG}_n(\vbeta_n)&=&-\frac{1}{2}\sum_{i=1}^n\sum_{j=1}^m\bigl(1-2\pi_{ij}(\vbeta_n)\bigr)\vA_{ij}^{1/2}(\vbeta_n)
\vX_{ij}\vX_{ij}^T\ve_j^T\overline{\vR}^{-1}\vepsilon_i(\vbeta_n),
\end{eqnarray*}
where
$\epsilon_{ij}(\vbeta_n)=\vA_{ij}^{-1/2}(\vbeta_n)
(Y_{ij}-\pi_{ij}(\vbeta_{n}))$,
$\vepsilon_i(\vbeta_n)=\vA_i^{-1/2}(\vbeta_n)(\vY_i-\vpi(\vbeta_n))$ and
$\ve_j$ denotes a unit vector of length $m$ whose $j$th entry is 1 and all
 other entries of which are 0.
 \end{lemma}

\begin{lemma} \label{tiger3}
Assume conditions \textup{(A1)}--\textup{(A4)}. If $n^{-1}p_n^2=o(1)$, then
$\forall \Delta>0$, for $\vb_n\in R^{p_n}$, we have
\begin{eqnarray*}
\sup_{\Vert\vbeta_n-\vbeta_{n0}\Vert\leq \Delta\sqrt{p_n/n}}\
\sup_{\Vert\vb_n\Vert=1}|\vb_n^T[\vD_n(\vbeta_n)-\overline{\vD}_n(\vbeta_{n})]\vb_n|
=O_p\bigl(\sqrt{np_n}\bigr).
\end{eqnarray*}
\end{lemma}

\begin{remark}\label{rem2}
The matrix $\vD_n(\vbeta_n)-\overline{\vD}_n(\vbeta_{n})$
is symmetric. The above lemma immediately implies that
\begin{eqnarray*}
\sup_{\Vert\vbeta_n-\vbeta_{n0}\Vert\leq \Delta\sqrt{p_n/n}}
|\lambda_{\min}[\vD_n(\vbeta_n)-\overline{\vD}_n(\vbeta_{n})]|&=&O_p\bigl(\sqrt{np_n}\bigr),\\
 \sup_{\Vert\vbeta_n-\vbeta_{n0}\Vert\leq
\Delta\sqrt{p_n/n}}|\lambda_{\max}[\vD_n(\vbeta_n)-\overline{\vD}_n(\vbeta_{n})]|&=&O_p\bigl(\sqrt{np_n}\bigr).
\end{eqnarray*}

Furthermore, we can use the leading term $\overline{\vH}_n(\vbeta_n)$ in
(\ref{rep}) to approximate the negative gradient function
$\overline{\vD}_n(\vbeta_n)$. This result is given by Lemma \ref{tiger4}
below. Lemma \ref{tiger5} further establishes an equicontinuity  result for
$\overline{\vH}_n(\vbeta_n)$.
\end{remark}

\begin{lemma} \label{tiger4}
Assume conditions \textup{(A1)}--\textup{(A4)}. If $n^{-1}p_n^2=o(1)$, then
$\forall \Delta>0$, for $\vb_n\in R^{p_n}$, we have
\[
\sup_{||\vbeta_n-\vbeta_{n0}\Vert\leq \Delta\sqrt{p_n/n}}\
\sup_{\Vert\vb_n\Vert=1}|\vb_n^T[\overline{\vD}_n(\vbeta_n)-\overline{\vH}_n(\vbeta_n)]\vb_n|
=O_p\bigl(\sqrt{n}p_n\bigr).
\]
\end{lemma}

\begin{lemma} \label{tiger5}
Assume conditions \textup{(A1)}--\textup{(A4)}. If $n^{-1}p_n^2=o(1)$, then
$\forall \Delta>0$, for $\vb_n\in R^{p_n}$, we have
\[
\sup_{\Vert\vbeta_n-\vbeta_{n0}\Vert\leq \Delta\sqrt{p_n/n}}\
\sup_{\Vert\vb_n\Vert=1}|\vb_n^T[\overline{\vH}_n(\vbeta_n)-\overline{\vH}_n(\vbeta_{n0})]\vb_n|
=O_p\bigl(\sqrt{n}p_n\bigr).
\]
\end{lemma}

The proofs of Lemmas \ref{tiger1}--\ref{tiger4} are given in the \hyperref[app]{Appendix}; the proof of Lemma~\ref{tiger5}
is given in the supplementary article [Wang (\citeyear{wan2010})]. The
following theorem ensures the existence and consistency of the GEE estimator
when $p_n\ra \infty$.

\begin{theorem}[(Existence and  consistency)]\label{thm1}
Assume conditions
\textup{(A1)}--\textup{(A4)} and that $n^{-1}p_n^2=o(1)$. Then, $\vS_n(\vbeta_n)=0$ has a root
$\widehat{\vbeta}_n$ such that
\[
\Vert\widehat{\vbeta}_n-\vbeta_{n0}\Vert=O_p\bigl(\sqrt{p_n/n}\bigr).
\]
\end{theorem}

\begin{pf}
We will prove that (\ref{bull}) holds.  This requires us to
evaluate the sign of $(\vbeta_n-\vbeta_{n0})^T\vS_n(\vbeta_n)$ on
$\{\vbeta_n\dvtx
\Vert\vbeta_n-\vbeta_{n0}\Vert= \Delta\sqrt{p_n/n}\}$. Note that
\begin{eqnarray*}
&&(\vbeta_n-\vbeta_{n0})^T\vS_n(\vbeta_n)\\
&&\qquad=(\vbeta_n-\vbeta_{n0})^T\vS_n(\vbeta_{n0})-(\vbeta_n-\vbeta_{n0})^T\vD_n(\vbeta_n^*)(\vbeta_n-\vbeta_{n0})\\
&&\qquad\triangleq  I_{n1}+I_{n2},
\end{eqnarray*}
where $\vbeta_n^*$ lies between $\vbeta_n$
and $\vbeta_{n0}$, that is, $\vbeta_n^*=t\vbeta_n+(1-t)\vbeta_{n0}$ for some
$0<t<1$. Next, we write
\begin{eqnarray*}
I_{n1}&=&(\vbeta_n-\vbeta_{n0})^T\overline{\vS}_n(\vbeta_{n0})+
(\vbeta_n-\vbeta_{n0})^T[\vS_n(\vbeta_{n0})-\overline{\vS}_n(\vbeta_{n0})]\\
&\triangleq& I_{n11}+I_{n12}.
\end{eqnarray*}
We have $|I_{n11}|\leq \Vert\vbeta_n-\vbeta_{n0}\Vert\cdot
\Vert\overline{\vS}_n(\vbeta_{n0})\Vert=\Delta\sqrt{p_n/n}\Vert\overline{\vS}_n(\vbeta_{n0})\Vert
$ by the Cauchy--Schwarz inequality.  Furthermore,
\begin{eqnarray*}
&&E[\Vert\overline{\vS}_n(\vbeta_{n0})\Vert^2]\\
&&\qquad=E\Biggl\{\sum_{i=1}^n\vepsilon_i^T(\vbeta_{n0})\overline{\vR}^{-1}\vA_{i}^{1/2}(\vbeta_{n0})
\vX_i\vX_i^T\vA_{i}^{1/2}(\vbeta_{n0})\overline{\vR}^{-1}\vepsilon_i(\vbeta_{n0})\Biggr\}\\
&&\qquad\leq\sum_{i=1}^n\lambda_{\max}(\vX_i\vX_i^T)\lambda_{\max}(\vA_{i}(\vbeta_{n0}))
\lambda_{\max}(\overline{\vR}^{-2})E[\vepsilon_i^T(\vbeta_{n0})\vepsilon_i(\vbeta_{n0})]\\
&&\qquad\leq  C\operatorname{Tr}\Biggl(\sum_{i=1}^n\vX_i\vX_i^T\Biggr)=C\sum_{i=1}^n\sum_{j=1}^m\vX_{ij}^T\vX_{ij}=O(np_n).
\end{eqnarray*}
Here, and throughout the paper, we use $C$ to denote a generic positive
constant which may vary from line to line. Thus,
$\Vert\overline{\vS}_n(\vbeta_{n0})\Vert=O_p(\sqrt{np_n})$. This implies that
$|I_{n11}|=\Delta O_p(p_n)$. For $I_{n12}$, we have
\[
|I_{n12}|\leq
\Vert\vbeta_n-\vbeta_{n0}\Vert\cdot
\Vert\vS_n(\vbeta_{n0})-\overline{\vS}_n(\vbeta_{n0})\Vert=\Delta
\sqrt{p_n/n}O_p(p_n)=\Delta o_p(p_n),
\]
by Lemma \ref{tiger1}. Hence,
$|I_{n1}|=\Delta O_p(p_n)$. In what follows, we evaluate $I_{n2}$:
\begin{eqnarray*}
I_{n2}&=&-(\vbeta_n-\vbeta_{n0})^T\overline{\vD}_n(\vbeta_n^*)(\vbeta_n-\vbeta_{n0})\\
&&-(\vbeta_n-\vbeta_{n0})^T[\vD_n(\vbeta_n^*)-\overline{\vD}_n(\vbeta_n^*)](\vbeta_n-\vbeta_{n0})\\
&\triangleq & I_{n21}+I_{n22}.
\end{eqnarray*}
First, note that
\begin{eqnarray*}
|I_{n22}|&\leq&\max\bigl(\bigl|\lambda_{\max}\bigl(\vD_n(\vbeta_n^*)-\overline{\vD}_n(\vbeta_n^*)\bigr)\bigr|,
\bigl|\lambda_{\min}\bigl(\vD_n(\vbeta_n^*)-\overline{\vD}_n(\vbeta_n^*)\bigr)\bigr|\bigr)\\
&&{}\times\Vert\vbeta_n-\vbeta_{n0}\Vert^2\\
&=&O_p\bigl(\sqrt{np_n}\bigr)\Delta^2 \frac{p_n}{n}= \Delta^2 o_p(p_n),
\end{eqnarray*}
by Lemma
\ref{tiger3}. On the other hand,
\begin{eqnarray*}
I_{n21}&=&-(\vbeta_n-\vbeta_{n0})^T\overline{\vH}_n(\vbeta_{n0})(\vbeta_n-\vbeta_{n0})\\
&&-(\vbeta_n-\vbeta_{n0})^T[\overline{\vH}_n(\vbeta_{n}^*)-\overline{\vH}_n(\vbeta_{n0})](\vbeta_n-\vbeta_{n0})\\
&&-(\vbeta_n-\vbeta_{n0})^T[\overline{\vD}_n(\vbeta_{n}^*)-\overline{\vH}_n(\vbeta_{n}^*)](\vbeta_n-\vbeta_{n0})\\
&\triangleq &  I_{n21}^a+ I_{n21}^b+ I_{n21}^c.
\end{eqnarray*}
From Lemma
\ref{tiger5}, we have $I_{n21}^b=\Delta^2 o_p(p_n)$;  from Lemma
\ref{tiger4}, we have $I_{n21}^c=\Delta^2 o_p(p_n)$. Finally, we evaluate $I_{n21}^a$.
We have
\begin{eqnarray*}
I_{n21}^a&=&-(\vbeta_n-\vbeta_{n0})^T
\Biggl[\sum_{i=1}^n\vX_i^T\vA_i^{1/2}(\vbeta_{n0})\overline{\vR}^{-1}\vA_i^{1/2}(\vbeta_{n0})\vX_i\Biggr]
(\vbeta_n-\vbeta_{n0})\\
&\leq &
-\lambda_{\min}(\overline{\vR}^{-1})\min_{i}\lambda_{\min}(\vA_i(\vbeta_{n0}))
\lambda_{\min}\Biggl(\sum_{i=1}^n\vX_i^T\vX_i\Biggr)\Vert\vbeta_n-\vbeta_{n0}\Vert^2\\
&\leq & -C \Delta^2 p_n,
\end{eqnarray*}
by (A3). Thus,
$(\vbeta_n-\vbeta_{n0})^T\vS_n(\vbeta_n)$ on $\{\vbeta_n\dvtx
\Vert\vbeta_n-\vbeta_{n0}\Vert= \Delta\sqrt{p_n/n}\}$ is asymptotically dominated in
probability by $I_{n11}+I_{n21}^a$, which is negative for $\Delta$ large
enough.
\end{pf}

\subsection{Asymptotic normality of the GEE estimator}\label{sec33}

The asymptotic distribution of the GEE estimator $\widehat{\vbeta}_n$ is
closely related to that of the ideal estimating function
$\overline{\vS}_n(\vbeta_{n0})$. When appropriately normalized,
$\overline{\vS}_n(\vbeta_{n0})$ has an asymptotic normal distribution, as
shown by the following lemma.

\begin{lemma}\label{tiger6}
Assume conditions \textup{(A1)}--\textup{(A4)}. If $n^{-1}p_n^3=o(1)$, then
$\forall \valpha_n\in R^{p_n}$ such that $\Vert\valpha_n\Vert=1$, we have
\begin{eqnarray*}
\valpha_n^T\overline{\vM}^{-1/2}_n(\vbeta_{n0})\overline{\vS}_n(\vbeta_{n0})\ra
N(0,1)\quad\mbox{in distribution}.
\end{eqnarray*}
\end{lemma}

To prove Lemma \ref{tiger6}, we write
$\valpha_n^T\overline{\vM}^{-1/2}_n(\vbeta_{n0})\overline{\vS}_n(\vbeta_{n0})$
as a sum of independent random variables and then check the Lindberg--Feller
condition for the central limit theorem. The detailed proof is given in the
\hyperref[app]{Appendix}. The following theorem ensures the asymptotic normality of the GEE
estimator when $n^{-1}p_n^3=o(1)$.

\begin{theorem}[(Asymptotic normality)]\label{anorm}
Assume conditions \textup{(A1)}--\textup{(A4)}. If
$n^{-1}p_n^3=o(1)$, then $\forall \valpha_n\in R^{p_n}$ such that
$\Vert\valpha_n\Vert=1$, we have
\begin{eqnarray*}
\valpha_n^T\overline{\vM}^{-1/2}_n(\vbeta_{n0})\overline{\vH}_n(\vbeta_{n0})(\widehat{\vbeta}_{n}-\vbeta_{n0})\ra
N(0,1)
\end{eqnarray*}
in distribution.
\end{theorem}

\begin{pf}
We have
\begin{eqnarray*}
&&\valpha_n^T\overline{\vM}^{-1/2}_n(\vbeta_{n0})\overline{\vS}_n(\vbeta_{n0})\\
&&\qquad= \valpha_n^T\overline{\vM}^{-1/2}_n(\vbeta_{n0})\vS_n(\vbeta_{n0})+
\valpha_n^T\overline{\vM}^{-1/2}_n(\vbeta_{n0})[\overline{\vS}_n(\vbeta_{n0})-\vS_n(\vbeta_{n0})]\\
&&\qquad=\valpha_n^T\overline{\vM}^{-1/2}_n(\vbeta_{n0})\vD_n(\vbeta_n^*)(\widehat{\vbeta}_{n}-\vbeta_{n0})\\
&&\qquad\quad{}+\valpha_n^T\overline{\vM}^{-1/2}_n(\vbeta_{n0})[\overline{\vS}_n(\vbeta_{n0})-\vS_n(\vbeta_{n0})]\\
&&\qquad=\valpha_n^T\overline{\vM}^{-1/2}_n(\vbeta_{n0})\overline{\vH}_n(\vbeta_{n0})(\widehat{\vbeta}_{n}-\vbeta_{n0})\\
&&\qquad\quad{}+\valpha_n^T\overline{\vM}^{-1/2}_n(\vbeta_{n0})[\vD_n(\vbeta_n^*)-\overline{\vH}_n(\vbeta_{n0})](\widehat{\vbeta}_{n}-\vbeta_{n0})\\
&&\qquad\quad{}+\valpha_n^T\overline{\vM}^{-1/2}_n(\vbeta_{n0})[\overline{\vS}_n(\vbeta_{n0})-\vS_n(\vbeta_{n0})],
\end{eqnarray*}
where, to obtain the second equality, we note that
$\vS_n(\widehat{\vbeta}_{n})=0$ and thus, by a Taylor expansion,
$\vS_n(\vbeta_{n0})=\vD_n(\vbeta_n^*)(\widehat{\vbeta}_{n}-\vbeta_{n0})$ for
some $\vbeta_n^*$ between $\widehat{\vbeta}_{n}$ and $\vbeta_{n0}$. By Lemma
\ref{tiger6}, $\valpha_n^T\overline{\vM}^{-1/2}_n(\vbeta_{n0})\overline{\vS}_n(\vbeta_{n0})\ra
N(0,1)$.   Therefore, to prove the theorem, it is sufficient to verify that
$\forall \Delta>0$,
\begin{eqnarray}\label{dog1}
&&\sup_{\Vert\vbeta_n-\vbeta_{n0}\Vert\leq
\Delta\sqrt{p_n/n}}|\valpha_n^T\overline{\vM}^{-1/2}_n(\vbeta_{n0})[\vD_n(\vbeta_n)-
\overline{\vH}_n(\vbeta_{n0})](\widehat{\vbeta}_{n}-\vbeta_{n0})|\nonumber\\
[-8pt]\\ [-8pt]
&&\qquad=o_p(1)\nonumber
\end{eqnarray}
and\vspace*{1pt}
\begin{equation}\label{dog2}
|\valpha_n^T\overline{\vM}^{-1/2}_n(\vbeta_{n0})
[\overline{\vS}_n(\vbeta_{n0})-\vS_n(\vbeta_{n0})]|=o_p(1).
\end{equation}\vspace*{1pt}
We prove
(\ref{dog2}) first. Note that
\begin{eqnarray*}
&&\bigl[\valpha_n^T\overline{\vM}^{-1/2}_n(\vbeta_{n0})
[\overline{\vS}_n(\vbeta_{n0})-\vS_n(\vbeta_{n0})]\bigr]^2\\
&&\qquad=\valpha_n^T\overline{\vM}^{-1/2}_n(\vbeta_{n0})
[\overline{\vS}_n(\vbeta_{n0})-\vS_n(\vbeta_{n0})]
[\overline{\vS}_n(\vbeta_{n0})-\vS_n(\vbeta_{n0})]^T\overline{\vM}^{-1/2}_n(\vbeta_{n0})\valpha_n\\
&&\qquad\leq \lambda_{\max}(\overline{\vM}^{-1}_n(\vbeta_{n0}))\lambda_{\max}\bigl([\overline{\vS}_n(\vbeta_{n0})-\vS_n(\vbeta_{n0})]
[\overline{\vS}_n(\vbeta_{n0})-\vS_n(\vbeta_{n0})]^T\bigr)\\
&&\qquad\leq\frac{\Vert\overline{\vS}_n(\vbeta_{n0})-\vS_n(\vbeta_{n0})\Vert^2}{\lambda_{\min}(\overline{\vM}_n(\vbeta_{n0}))}\\
&&\qquad\leq\frac{\Vert\overline{\vS}_n(\vbeta_{n0})-\vS_n(\vbeta_{n0})\Vert^2}{C\lambda_{\min}(\sum_{i=1}^n\vX_i^T\vX_i)}\\
&&\qquad=O_p(p_n^2/n)=o_p(1),
\end{eqnarray*}
by Lemma \ref{tiger1} and the fact that
%
\begin{equation}\label{mn}
\lambda_{\min}(\overline{\vM}_n(\vbeta_{n0}))\geq
C\lambda_{\min}\Biggl(\sum_{i=1}^n\vX_i^T\vX_i\Biggr).
\end{equation}
A justification of (\ref{mn})
is given in the proof of Lemma \ref{tiger6} in the \hyperref[app]{Appendix}. Thus, (\ref{dog2})
holds. Next,\vspace*{1pt} we prove (\ref{dog1}).  We have
\begin{eqnarray*}
&&\sup_{\Vert\vbeta_n-\vbeta_{n0}\Vert\leq
\Delta\sqrt{p_n/n}}|\valpha_n^T\overline{\vM}^{-1/2}_n(\vbeta_{n0})[\vD_n(\vbeta_n)-
\overline{\vH}_n(\vbeta_{n0})](\widehat{\vbeta}_{n}-\vbeta_{n0})|\\
&&\qquad\leq \sup_{\Vert\vbeta_n-\vbeta_{n0}\Vert\leq
\Delta\sqrt{p_n/n}}|\valpha_n^T\overline{\vM}^{-1/2}_n(\vbeta_{n0})[\vD_n(\vbeta_n)-
\overline{\vD}_n(\vbeta_{n})](\widehat{\vbeta}_{n}-\vbeta_{n0})|\\
&&\qquad\quad{}+\sup_{\Vert\vbeta_n-\vbeta_{n0}\Vert\leq
\Delta\sqrt{p_n/n}}|\valpha_n^T\overline{\vM}^{-1/2}_n(\vbeta_{n0})[\overline{\vD}_n(\vbeta_n)-
\overline{\vH}_n(\vbeta_{n})](\widehat{\vbeta}_{n}-\vbeta_{n0})| \\
&&\qquad\quad{}+\sup_{\Vert\vbeta_n-\vbeta_{n0}\Vert\leq
\Delta\sqrt{p_n/n}}|\valpha_n^T\overline{\vM}^{-1/2}_n(\vbeta_{n0})[\overline{\vH}_n(\vbeta_n)-
\overline{\vH}_n(\vbeta_{n0})](\widehat{\vbeta}_{n}-\vbeta_{n0})| \\
&&\qquad\triangleq  I_{n1}+I_{n2}+I_{n3}.
\end{eqnarray*}\vspace*{1pt}
By the Cauchy--Schwarz inequality and
Remark \ref{rem2}, we have
\begin{eqnarray*}
I_{n1}&\leq&
\sup_{\Vert\vbeta_n-\vbeta_{n0}\Vert\leq
\Delta\sqrt{p_n/n}}\bigl[\valpha_n^T\overline{\vM}^{-1/2}_n(\vbeta_{n0})\bigl(\vD_n(\vbeta_n)-
\overline{\vD}_n(\vbeta_{n})\bigr)^2\\
&&\hspace*{161pt}{}\times\overline{\vM}^{-1/2}_n(\vbeta_{n0})\valpha_n\bigr]^{1/2}\Vert\widehat{\vbeta}_{n}-\vbeta_{n0}\Vert\\
&\leq & \sup_{\Vert\vbeta_n-\vbeta_{n0}\Vert\leq
\Delta\sqrt{p_n/n}}\max\bigl(\bigl|\lambda_{\min}(\vD_n(\vbeta_n)-
\overline{\vD}_n(\vbeta_{n})\bigr)\bigr|,\\
&&\hspace*{99pt}\bigl|\lambda_{\max}\bigl(\vD_n(\vbeta_n)-
\overline{\vD}_n(\vbeta_{n})\bigr)\bigr|\bigr)\\
&&\hphantom{\sup_{\Vert\vbeta_n-\vbeta_{n0}\Vert\leq
\Delta\sqrt{p_n/n}}}{}\times\lambda_{\min}^{-1/2}(\overline{\vM}_n(\vbeta_{n0}))
O_p(p_n^{1/2}n^{-1/2})\\
&=&
O_p\bigl(\sqrt{n}p_n\bigr)O(n^{-1/2})O_p(n^{-1/2}p_n^{1/2})=O_P(n^{-1/2}p_n^{3/2})=o_p(1).
\end{eqnarray*}
Hence, $I_{n1}=o_p(1)$. By the same argument and Lemmas
\ref{tiger4} and \ref{tiger5}, we
also have $I_{n2}=o_p(1)$ and $I_{n3}=o_p(1)$. This proves (\ref{dog1}).
\end{pf}

\begin{remark}
Note that the condition $n^{-1}p_n^3=o(1)$ is the same as
that of Huber~(\citeyear{hub1973}) for an M-estimator with independent data and diverging
number of parameters. It is weaker than the condition $n^{-1}p_n^5=o(1)$ in
Fan and Peng (\citeyear{fanpen2004}) and Lam and Fan (\citeyear{lamfan2008}) for asymptotic normality.
\end{remark}

\begin{remark}
Combining the result of Theorem \ref{anorm} with the
Cram\'er--Wold device, it is easy to see that for any $l\times p_n$ matrix
$\vB_n$ with $l$ fixed and such that $\vB_n\vB_n^T\ra \vF$, a
positive definite matrix, we have
\[
\vB_n\vSigma_n^{-1/2}(\vbeta_{n0})(\widehat{\vbeta}_n-\vbeta_{n0})\ra
\vN_{l}(\vnull, \vF),
\]
where
\[
\vSigma_n=\overline{\vH}^{-1}_n(\vbeta_{n0})\overline{\vM}_n(\vbeta_{n0})\overline{\vH}^{-1}_n(\vbeta_{n0}).
\]
Now, take $\vB_n=(\vL_n\vSigma_n\vL_n^T)^{-1/2}\vL_n\vSigma_n^{1/2}$,
where $\vL_n$ is an $l\times p_n$ matrix such that $\vL_n\vSigma_n\vL_n^T$ is
invertible. Then, $\vB_n\vB_n^T=\vI_l$ and we have the following corollary
which gives the asymptotic distribution of
$\vL_n(\widehat{\vbeta}_n-\vbeta_{n0})$.
\end{remark}

\begin{corollary}\label{cor39}
Under the same conditions as in
Theorem \ref{anorm}, if $n^{-1}p_n^3=o(1)$, then
\[
(\vL_n\vSigma_n\vL_n^T)^{-1/2}\vL_n(\widehat{\vbeta}_n-\vbeta_{n0})\ra
\vN_{l}(\vnull, \vI_l)
\]
in distribution.
\end{corollary}

\subsection{Sandwich covariance formula and large-sample Wald test}

Theorem \ref{anorm} and Corollary \ref{cor39} suggest that the covariance matrix of
$\widehat{\vbeta}_n$ is approximately $\vSigma_n$. To estimate $\vSigma_n$,
Liang and Zeger (\citeyear{liazeg1986}) proposed, in the ``fixed $p$'' setup,
the following well-known \textit{sandwich covariance matrix estimator:}
\[
\widehat{\vSigma}_n=\vH^{-1}_n(\widehat{\vbeta}_{n})
\widehat{\vM}_n(\widehat{\vbeta}_{n})\vH^{-1}_n(\widehat{\vbeta}_{n}),
\]
where $\vH_n(\vbeta_n)$ is defined similarly as $\overline{\vH}_n(\vbeta_n),$
but with $\overline{\vR}$ replaced by $\widehat{\vR}$;
$\widehat{\vM}_n(\vbeta_n)$ is defined similarly as
$\overline{\vM}_n(\vbeta_n),$ except that $\overline{\vR}$ is replaced by
$\widehat{\vR}$ and the unknown true correlation matrix $\vR_0$ is replaced by
$\vepsilon_i(\vbeta_n)\vepsilon_i^T(\vbeta_n)$, with $\vepsilon_i(\vbeta_n)$
defined in Lemma \ref{tiger2}. Based on Corollary \ref{cor39} and the sandwich
covariance matrix estimator, an asymptotic $(1-\alpha)$\% confidence interval
($0<\alpha<1$) for $\beta_j$ is
\begin{equation}\label{conf}
\widehat{\vbeta}_{j}\pm
z_{\alpha/2}\vu_j^T\widehat{\vSigma}_n\vu_j,
\end{equation}
where
$z_{\alpha/2}$ denotes the upper $\frac{\alpha}{2}$ quantile of the
standard normal distribution and $\vu_j$ is the unit vector of length $p_n$
with the $j$th element equal to 1 and all the other elements equal to 0.

The sandwich covariance formula plays an important role in GEE methodology. In
the ``fixed $p$'' setup, it is known that the sandwich covariance matrix
estimator provides a consistent estimator for the variance of the GEE
estimator, even when the working correlation matrix is misspecified. The
following theorem shows that this appealing property is still valid when $p_n$
converges to $\infty$ at an appropriate rate.

The proofs of Theorem \ref{cat5} and Corollary \ref{wald} below are given in
the \hyperref[app]{Appendix}.

\begin{theorem}\label{cat5}
Assume conditions \textup{(A1)}--\textup{(A4)} and that
$n^{-1}p_n^3=o(1)$. Then,
\[
\vC_n\widehat{\vSigma}_n\vC_n^T-\vC_n\vSigma_n\vC_n^T=o_p(n^{-1}),
\]
where
$\vC_n$ is any $l\times q_n$ matrix such that $l$ is fixed and
$\vC_n\vC_n^T=\vG$ with $\vG$ being an $l\times l$ positive definite matrix.
\end{theorem}

\begin{remark}
It is worth pointing out a subtle issue that is sometimes
overlooked in the existing literature on high-dimensional analysis of
independent data. In order to justify the validity of the asymptotic
confidence interval or large-sample test for estimable contrast, it is
necessary to show that the convergence rate in Theorem \ref{cat5} is
$o_p(n^{-1})$. Note that the estimable contrast is asymptotically normal with
convergence rate $O_p(n^{1/2})$; see, for example, Corollary 2.1 in He and
Shao~(\citeyear{hesha2000}) for the case of an M-estimator based on independent data. In the
literature, sometimes only the $o_p(1)$ rate is provided, which is not adequate,
but can be fixed.

Next, we consider the large-sample Wald test for testing the following linear
hypothesis:
\[
H_0\dvtx \vL_n\vbeta_{n0}=\vnull\quad\mbox{vs.}\quad H_1\dvtx\vL_n\vbeta_{n0}\neq \vnull,
\]
where $\vL_n$ is an $l\times p_n$ matrix
with $l$ fixed and $\vL_n\vL_n^T=\vI_l$. The Wald test statistic is defined as
\[
W_n=(\vL_n\widehat{\vbeta}_n)^T(\vL_n\widehat{\vSigma}_n\vL_n^T)^{-1}(\vL_n\widehat{\vbeta}_n).
\]
The corollary below shows that the Wald test remains valid, even when the
number of covariates diverges with the sample size.
\end{remark}

\begin{corollary}\label{wald}
Assume conditions \textup{(A1)}--\textup{(A4)}. If $n^{-1}p_n^3=o(1)$, then
$W_n\ra \chi^2_l$ in distribution under $H_0$, where $\chi^2_l$ denotes the
$\chi^2$ distribution with $l$ degrees of freedom.
\end{corollary}

\begin{remark}\label{rem6}
For testing a high-dimensional hypothesis $H_0\dvtx
\vbeta_n=\vbeta_{n0}^*$ versus $H_1\dvtx \vbeta_n\neq\vbeta_{n0}^*$, it can be
shown that
\begin{eqnarray}\label{high}
\frac{(\widehat{\vbeta}_n-\vbeta_{n0}^*)^T\widehat{\vSigma}_n^{-1}(\widehat{\vbeta}_n-\vbeta_{n0}^*)-p_n}{\sqrt{2p_n}}\ra
N(0,1)
\end{eqnarray}
in distribution under $H_0,$ under some regularity conditions. A
proof of this result is given in the supplementary article [Wang (\citeyear{wan2010})].
\end{remark}

\section{Numerical studies}\label{sec4}
We consider the following model for the marginal expectation of $Y_{ij}$,
$i=1,\ldots,n$, given $\vX_{ij}$,
\begin{equation}\label{simu}
\operatorname{logit}(\pi_{ij})=X_{ij}^T\vbeta_{n0},\quad j=1, 2, 3,
\end{equation}
where
$\vbeta_{n0}$ is a $p_n$-dimensional vector of parameters with
$p_n=\lfloor 2.5n^{1/3}\rfloor$, with $\lfloor q
\rfloor$ denoting the the largest integer not greater than $q$. In this
example, $\vbeta_{n0}^T=(0.4\cdot\veins_k^T,-0.3\cdot \veins_k^T, 0.2\cdot
\veins_k^T, -0.1\cdot \veins_{p_n-3k}^T)$, where $\veins_k$ denotes a
$k$-dimensional vector of 1's and $k=\lfloor p_n/4 \rfloor$. In
addition, $X_{ij}=(x_{ij1},\ldots,x_{ijp_n})^T$ has a multivariate normal
distribution with mean zero, marginal variance 0.2 and an AR-1 correlation
matrix with autocorrelation coefficient 0.5. The binary response vector for
each cluster has the above marginal mean and an exchangeable (also called
compound symmetry or CS) correlation structure with correlation coefficient
0.5. Such correlated binary data are generated using Bahadur's representation
[see, e.g., Fitzmaurice (\citeyear{fit1995})].

Since, for different sample sizes, the parameter dimension is different, we
measure the accuracy of estimation by the \textit{simulated average mean square
error}, which is obtained by averaging
$\Vert\widehat{\vbeta}_n-\vbeta_{n0}\Vert^2/p_n$ over 500 simulated samples. Table
\ref{tab1}
reports simulation results using four different working correlation
structures: independence working correlation matrix (IN), unstructured working
correlation matrix (UN), compound symmetry working correlation matrix (CS) and
the first order autocorrelation working correlation matrix (AR-1), for sample
sizes $n=500$, $1000$, $2000$ and $3000$. Table \ref{tab1} demonstrates that when the
covariate dimension grows at an appropriate rate with the sample size, the
accuracy of GEE estimator is satisfactory. We also observe that when the true
correlation matrix (CS in this case) is adopted, the estimator is more
efficient.

\begin{table}\tablewidth=260pt
\caption{The simulated average mean squared error ($\times$10) for estimating
$\vbeta_{n0}$ using four different working correlation structures}\label{tab1}
\begin{tabular*}{260pt}{@{\extracolsep{\fill}}lccccc@{}}
\hline
& & \multicolumn{4}{c@{}} {\textbf{Working correlation structure}} \\
[-7pt]
& &\multicolumn{4}{c@{}}{\hrulefill}\\
$\bolds n$ & $\bolds{p_n}$ & \textbf{IN} &
\textbf{UN} & \textbf{CS} & \textbf{AR-1} \\
\hline
\phantom{0}500 & 19 & 0.265 & 0.156 & 0.154 & 0.179 \\
1000 & 24 & 0.141 & 0.103 & 0.100 & 0.111 \\
2000 & 31 & 0.090 & 0.074 & 0.071 & 0.075\\
3000 & 36 & 0.070 & 0.065 & 0.063 & 0.065\\
\hline
\end{tabular*}
\end{table}

We next examine the accuracy of the sandwich variance formula. The standard
deviations of the estimated coefficients over 500 simulations are averaged and
regarded as the true standard error (SD). Table \ref{sphericcase} compares SD with the
standard error obtained from the sandwich variance formula (SD2) when the
unstructured working correlation matrix is used for estimating
$\widehat{\beta}_k$, $\widehat{\beta}_{2k}$, $\widehat{\beta}_{3k}$ and
$\widehat{\beta}_{p_n}$. We observe that the sandwich variance formula works
remarkably well. Similar phenomena are also observed for estimating other
regression coefficients and with other working correlation structures, but, for reasons of brevity, these are
not reported.
\begin{table}[b]\tablewidth=300pt
\tabcolsep=0pt
\caption{Standard deviation (SD) and estimated standard deviation (SD2) using
the sandwich variance formula} \label{sphericcase}
\begin{tabular*}{300pt}{@{\extracolsep{\fill}}lccccccccc@{}}
\hline
& & \multicolumn{2}{c}{$\bolds{\widehat{\beta}_k}$} & \multicolumn{2}{c} {$\bolds{\widehat{\beta}_{2k}}$}
& \multicolumn{2}{c} {$\bolds{\widehat{\beta}_{3k}}$} & \multicolumn{2}{c@{}}
{$\bolds{\widehat{\beta}_{p_n}}$}\\ [-7pt]
&&\multicolumn{2}{c}{\hrulefill} &\multicolumn{2}{c}{\hrulefill} &\multicolumn{2}{c}{\hrulefill} &\multicolumn{2}{c@{}}{\hrulefill}\\
$\bolds n$ &$\bolds{p_n}$ & \textbf{SD} &
\textbf{SD2} & \textbf{SD} & \textbf{SD2} & \textbf{SD}
& \textbf{SD2} & \textbf{SD} & \textbf{SD2}\\
\hline
\phantom{0}500 & 19 & 0.126 & 0.111 & 0.114 & 0.110 & 0.117 & 0.111  & 0.089 & 0.098\\
1000 & 24 & 0.082 & 0.083 & 0.079 & 0.083 & 0.085 & 0.083  &0.072 & 0.074\\
2000 & 31 & 0.073 & 0.060 & 0.063 & 0.060 & 0.065 & 0.060  & 0.051 & 0.053\\
3000 & 36 & 0.060 & 0.051 & 0.049 & 0.051 & 0.052 & 0.051  & 0.051 & 0.045\\
\hline
\end{tabular*}
\end{table}

\begin{figure}

\includegraphics{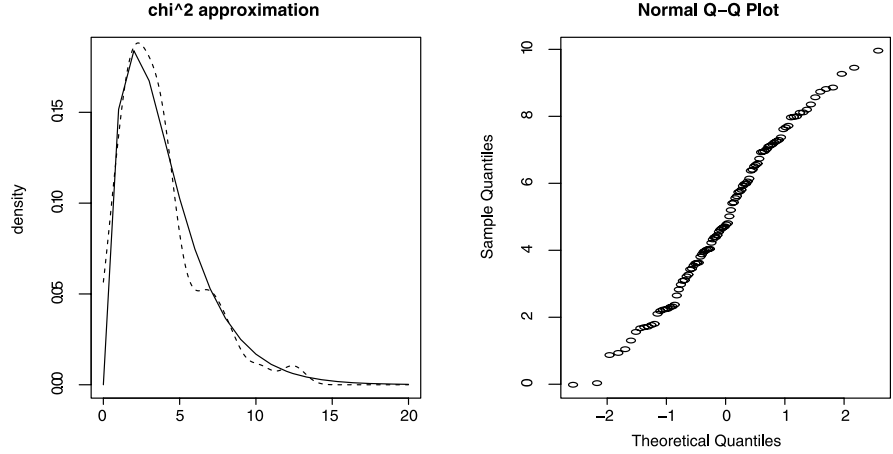}

\caption{The left panel gives the estimated null density of the
large-sample Wald test (dashed curve) and the density of the chi-square
distribution with four degrees of freedom (solid curve) for testing
$H_0\dvtx
\beta_{21}=\beta_{22}=\beta_{23}=\beta_{24}=0$. The right panel gives the
normal Q--Q plot of the Wald test statistic under the null hypothesis
$\vbeta_n=\vbeta_{n0}$.}\label{fig1}
\end{figure}

Finally, we investigate hypothesis testing based on the large-sample Wald
test. We consider model (\ref{simu}) with $n=1000$, $p_n=24$ and
$\vbeta_{n0}^T=(0.4\cdot\veins_6^T,-0.3\cdot \veins_6^T, 0.2\cdot \veins_6^T,
-0.1\cdot \veins_{2}^T, 0,0,0,0)$. The left panel of Figure \ref{fig1} depicts the
density of the Wald test under the null hypothesis $H_0\dvtx
\beta_{21}=\beta_{22}=\beta_{23}=\beta_{24}=0$ and compares it with the
density curve of the $\chi^2_4$ distribution. It demonstrates that the
$\chi^2$ approximation given in Corollary \ref{wald} is accurate. The right panel of
Figure~\ref{fig1} gives the normal Q--Q plot for the Wald test statistic under the null
hypothesis $\vbeta_n=\vbeta_{n0}$ and it shows that the null distribution can
be approximated well by a normal distribution for testing a higher-dimensional alternative, as discussed in Remark \ref{rem6}.

\section{Discussions}\label{sec5}
\subsection{Extension to general GEE}\label{sec51}
 Although the focus of the paper is on clustered binary data, the approaches
 and techniques can be extended to general GEE. For general GEE, the decomposition of $\overline{\vD}_n(\vbeta_n)$
given in Lemma \ref{tiger2} has a more complex expression, and the potential
unboundedness of $Y_{ij}$ makes the derivation of various probability bounds
and  asymptotic equivalence more delicate. Below, we give a brief discussion of
the large-$p$ asymptotics for general GEE.

Assume that the first two marginal moments of $Y_{ij}$ are
$\mu_{ij}(\vbeta_n):=\E_{\vbeta_n}(Y_{ij})=\mu(\theta_{ij})$ and
$\sigma^2_{ij}(\vbeta_n):=\Var_{\vbeta_n}(Y_{ij})=\dot{\mu}(\theta_{ij})$,
where $\theta_{ij}=\vX_{ij}^T\vbeta_n$. These moment assumptions would follow
when the marginal response variable has a canonical exponential family
distribution with scaling parameter 1. Let
$\vA_{i}(\vbeta_n)=\diag(\sigma^2_{i1}(\vbeta_n),\ldots,\sigma^2_{im}(\vbeta_n))$
and $\vmu_i(\vbeta_n)=(\mu_{i1}(\vbeta_n),\ldots,\mu_{im}(\vbeta_n))^T$. The
\textit{GEE estimator} $\widehat{\vbeta}_n$ is the solution of
\begin{equation}\label{ggee}
\sum_{i=1}^n\vX_i^T\vA_i^{1/2}(\vbeta_n)\widehat{\vR}^{-1}\vA_i^{-1/2}(\vbeta_n)\bigl(\vY_i-\vmu_i(\vbeta_n)\bigr)=0.
\end{equation}

In addition to assumptions (A1)--(A4) in Section \ref{sec32}, we adopt two additional
conditions:
\begin{longlist}
\item[(A5)] there exists a finite constant $M_1>0$ such that
$E(\Vert\vA_i^{-1/2}(\vbeta_n)(\vY_i-\vmu_i(\vbeta_n))\Vert^{2+\delta})\leq M_1$ for all $i$ and some $\delta>0$;\\

\item[(A6)] if  $B_n=\{\vbeta_n\dvtx \Vert\vbeta_n-\vbeta_{n0}\Vert\leq \Delta \sqrt{p_n/n}\}$,
then $\dot{\mu}(\vX_{ij}^T\vbeta_n)$, $1\leq i \leq n$, $1\leq j\leq m$, are
uniformly bounded away from 0 and $\infty$ on $B_n$;
$\ddot{\mu}(\vX_{ij}^T\vbeta_n)$ and
$\mu^{(3)}(\vX_{ij}^T\vbeta_n)$, $1\leq i \leq n$, $1\leq j\leq m$, are uniformly bounded by a finite positive constant $M_2$ on $B_n$.
\end{longlist}

\begin{remark}
Condition (A5) is similar to the condition in Lemma 2 of
Xie and Yang (\citeyear{xieyan2003}) and condition $(\tilde{N}_{\delta})$ in Balan and
Schiopu-Kratina (\citeyear{balsch2005}). Condition~(A6) requires
$\mu_{ij}^{(k)}(\vX_{ij}^T\vbeta_n)$, $k=1,2,3$, to be uniformly bounded when
$\vbeta_n$ is in a local neighborhood around $\vbeta_{n0}$. This condition is
generally satisfied for GEE. For example, when the marginal model follows a
Poisson distribution, $\mu(t)=\exp(t)$, thus
$\mu_{ij}^{(k)}(\vX_{ij}^T\vbeta_n)=\exp(\vX_{ij}^T\vbeta_n)$, $k=1,2,3$, are
uniformly bounded on  $B_n$.
\end{remark}

\begin{theorem}\label{general} Assume conditions \textup{(A1)}--\textup{(A6)} and that $n^{-1}p_n^2=o(1)$.
The generalized estimating equation (\ref{ggee}) then has a root
$\widehat{\vbeta}_n$ such that $
\Vert\widehat{\vbeta}_n-\vbeta_{n0}\Vert=O_p(\sqrt{p_n/n})$.  Furthermore, if
$n^{-1}p_n^3=o(1)$, then $\forall \valpha_n\in R^{p_n}$ such that
$\Vert\valpha_n\Vert=1$,
\[
\valpha_n^T\overline{\vM}^{-1/2}_n(\vbeta_{n0})\overline{\vH}_n(\vbeta_{n0})(\widehat{\vbeta}_{n}-\vbeta_{n0})\ra
N(0,1)
\]
in distribution, where $\overline{\vM}^{-1/2}_n(\vbeta_{n0})$ and
$\overline{\vH}_n(\vbeta_{n0})$ have the same expressions as in Section 3.2.
\end{theorem}

A sketch of the proof of Theorem \ref{general} is given in the supplementary article [Wang (\citeyear{wan2010})].

\subsection{Related problems}

In some scenarios, a ``large $n$, diverging $m$'' asymptotic framework, where
$p$ is either fixed or also diverges at an appropriate rate, may be more
appropriate. This corresponds to a real situation where the cluster size
is itself large. For example, in a longitudinal study, doctors take
measurements on the patients during each visit. Each patient forms a cluster.
The cluster size is large if the number of visits is large.
For a fixed $p$ setting, this ``large $n$, diverging $m$'' asymptotic framework
has been considered by Xie and Yang (\citeyear{xieyan2003}). A future topic of
interest is to consider large $m$ together with large $p$.

Another interesting direction for future study is to consider a more flexible
semiparametric specification for the generalized estimating equations in the
large-$p$ setting. In the classical ``fixed $p$'' setting, GEE with partially
linear model specification has been investigated by Lin and Carroll (\citeyear{lincar2001a}, \citeyear{lincar2001b}), Lin and
Ying (\citeyear{linyin2001}), He, Zhu and Fung (\citeyear{hezhufun2002}), Fan and Li (\citeyear{fanli2004}), Chiou and M\"uller
(\citeyear{chimul2005}), Wang, Carroll and Lin (\citeyear{wancarlin2005}), Chen and Jin (\citeyear{chejin2006}), He,
Fung and Zhu (\citeyear{hefunzhu2006}) and Huang, Zhang and Zhou (\citeyear{huazhazho2007}), among others.

\begin{appendix}
\section*{Appendix}\label{app}

We use $C$ to denote a generic positive constant that
can vary
from line to line.
\begin{pf*}{Proof of (\ref{inia})} It suffices [Ortega and
Rheinboldt (\citeyear{ortrhe1970})] to show\break that $\forall \epsilon>0$, there exists a
$\Delta>0$ such that for all $n$ sufficiently large,\break $
P(\sup_{\Vert\vbeta_n-\vbeta_{n0}\Vert=\Delta\sqrt{p_n/n}}(\vbeta_n-\vbeta_{n0})^T
\widetilde{S}_n(\vbeta_n)<0)\geq 1-\epsilon$. We
have
\begin{eqnarray*}
&&(\vbeta_n-\vbeta_{n0})^T
\widetilde{S}_n(\vbeta_n)\\
&&\qquad=(\vbeta_n-\vbeta_{n0})^T
\widetilde{S}_n(\vbeta_{n0})+(\vbeta_n-\vbeta_{n0})^T \frac{\partial}{\partial
\vbeta_{n}^T}\widetilde{S}_n(\vbeta_n^*)(\vbeta_n-\vbeta_{n0})\\
&&\qquad\triangleq  I_{n1}+I_{n2},
\end{eqnarray*}
where $\vbeta_n^*$ lies between
$\vbeta_{n0}$ and $\vbeta_n$. We first consider $I_{n1}$. For any $\vbeta_{n}$
such that $\Vert\vbeta_n-\vbeta_{n0}\Vert=\Delta\sqrt{\frac{p_n}{n}}$, we have
$|I_{n1}|\leq \Delta\sqrt{\frac{p_n}{n}}\Vert\widetilde{S}_n(\vbeta_{n0})\Vert$.
Note that
\begin{eqnarray*}
E[\Vert\widetilde{S}_n(\vbeta_{n0})\Vert^2]
&=&E\Biggl[\sum_{i=1}^n\bigl(\vY_i-\vpi_i(\vbeta_{n0})\bigr)^T\vX_i\vX_i^T\bigl(\vY_i-\vpi_i(\vbeta_{n0})\bigr)\Biggr]\\
&\leq & E\Biggl[\sum_{i=1}^n\lambda_{\max}(\vX_i\vX_i^T)\Vert\vY_i-\vpi_i(\vbeta_{n0})\Vert^2\Biggr]\\
&\leq& C \operatorname{Tr}\Biggl(\sum_{i=1}^n\vX_i\vX_i^T\Biggr)=C\sum_{i=1}^n\sum_{j=1}^m\vX_{ij}^T\vX_{ij}
=O(np_n),
\end{eqnarray*}
by assumption (A1). Thus, $|I_{n1}|\leq \Delta O_p(p_n)$. Next,
\begin{eqnarray*}
I_{n2}&=&-(\vbeta_n-\vbeta_{n0})^T\Biggl[\sum_{i=1}^n\vX_i^TA_i(\vbeta_{n0})\vX_i\Biggr](\vbeta_n-\vbeta_{n0})\\
&&-(\vbeta_n-\vbeta_{n0})^T\Biggl[\sum_{i=1}^n\vX_i^T\bigl(A_i(\vbeta^*)-A_i(\vbeta_{n0})\bigr)\vX_i\Biggr](\vbeta_n-\vbeta_{n0})\\
&\triangleq & I_{n21}+I_{n22}.
\end{eqnarray*}
Note that $I_{n21}\leq
-\lambda_{\min}(\vA_i(\vbeta_0))\lambda_{\min}(\sum_{i=1}^n\vX_i^T\vX_i)\Vert\vbeta_{n}-\vbeta_{n0}\Vert^2
\leq -C p_n \Delta^2,$ by~(A3). Since $\frac{\partial}{\partial
\vbeta_{n}}\vA_{ij}(\vbeta_n)=\pi_{ij}(\vbeta_n)(1-\pi_{ij}(\vbeta_n))(1-2\pi_{ij}(\vbeta_n))\vX_{ij}$, we have
\begin{eqnarray*}
|I_{n22}| & \leq &
(\vbeta_n-\vbeta_{n0})^T\Biggl[\sum_{i=1}^n\sum_{j=1}^m
|\vA_{ij}(\vbeta^*)-\vA_{ij}(\vbeta_{n0})|\vX_{ij}\vX_{ij}^T\Biggr](\vbeta_n-\vbeta_{n0})\\
&\leq & \sup_{i,j}\Vert\vX_{ij}\Vert\cdot \Vert\vbeta^*-\vbeta_0\Vert \cdot
\Vert\vbeta_n-\vbeta_0\Vert^2 \cdot\lambda_{\max}\Biggl(\sum_{i=1}^n\vX_i^T\vX_i\Biggr)\\
&\leq & O\bigl(\sqrt{p_n}\bigr)O_p\bigl(\sqrt{p_n/n}\bigr)(\Delta^2 p_n/n)O(n)=\Delta^2
o_p(p_n),
\end{eqnarray*}
by (A1)--(A3). Thus, for sufficiently
large $\Delta$, $(\vbeta_n-\vbeta_{n0})^T \widetilde{S}_n(\vbeta_n)$
is dominated by~$I_{n21}$, which is large and negative for all sufficiently
large $n$.
\end{pf*}

\begin{pf*}{Proof of (\ref{inib})}
The proof is given in the online supplement.
\end{pf*}

\begin{pf*}{Proof of Lemma \ref{tiger1}}
Let $\vQ=\{q_{j_1,j_2}\}_{1\leq j_1,j_2\leq m}$ denote the matrix
$\widehat{\vR}^{-1}-\overline{\vR}^{-1}$. Then,
\begin{eqnarray*}
&&\vS_n(\vbeta_{n0})-\overline{\vS}_n(\vbeta_{n0})\\
&&\qquad=\sum_{i=1}^n\sum_{j_1=1}^m\sum_{j_2=1}^m
q_{j_1,j_2}\vA^{1/2}_{ij_1}(\vbeta_{n0})\vA^{-1/2}_{ij_2}(\vbeta_{n0})\bigl(Y_{ij_2}-\pi_{ij_2}(\vbeta_{n0})\bigr)\vX_{ij_1}\\
&&\qquad= \sum_{j_1=1}^m\sum_{j_2=1}^mq_{j_1,j_2}
\Biggl[\sum_{i=1}^n\vA^{1/2}_{ij_1}(\vbeta_{n0})\epsilon_{ij_2}(\vbeta_{n0})\vX_{ij_1}\Biggr],
\end{eqnarray*}
where
$\epsilon_{ij_2}(\vbeta_{n0})=\vA^{-1/2}_{ij_2}(\vbeta_{n0})(Y_{ij_2}-\pi_{ij_2}(\vbeta_{n0}))$.
Note that
\begin{eqnarray*}
E\Biggl[\Biggl\Vert\sum_{i=1}^n\vA^{1/2}_{ij_1}(\vbeta_{n0})\epsilon_{ij_2}(\vbeta_{n0})\vX_{ij_1}\Biggr\Vert^2\Biggr]
&=&\sum_{i=1}^n\vA_{ij_1}(\vbeta_{n0})E[\epsilon^2_{ij_2}(\vbeta_{n0})]\vX_{ij_1}^T\vX_{ij_1}\\
&\leq& \sum_{i=1}^n\vX_{ij_1}^T\vX_{ij_1}=O(np_n).
\end{eqnarray*}
Thus, $ \Vert
\sum_{i=1}^n\vA^{1/2}_{ij_1}(\vbeta_{n0})\epsilon_{ij_2}(\vbeta_{n0})\vX_{ij_1}\Vert
=O_p(\sqrt{np_n})$ $\forall   1\leq j_1, j_2\leq m$ . Since, by (A4),
$q_{j_1,j_2}=O_p(\sqrt{p_n/n})$ $\forall  1\leq j_1, j_2\leq m $, the proof is complete.
\end{pf*}

\begin{pf*}{Proof of Lemma \ref{tiger2}}  The derivation can
be found in
Pan (\citeyear{pan2002}).
\end{pf*}

\begin{pf*}{Proof of Lemma \ref{tiger3}} Let $\vH_n(\vbeta_n)$,
$\vE_n(\vbeta_n)$ and $\vG_n(\vbeta_n)$ be defined the same as
$\overline{\vH}_n(\vbeta_n)$, $\overline{\vE}_n(\vbeta_n)$ and
$\overline{\vG}_n(\vbeta_n)$, respectively, but with $\overline{\vR}$ replaced
by $\widehat{\vR}$. By Lemma \ref{tiger2}, it is sufficient to prove the
\setcounter{equation}{0}
following three results:
\begin{eqnarray}\label{cat1}
&&\sup_{\Vert\vbeta_n-\vbeta_{n0}\Vert\leq\Delta\sqrt{p_n/n}}
\sup_{\Vert\vb_n\Vert=1}|\vb_n^T[\vH_n(\vbeta_n)-\overline{\vH}_n(\vbeta_{n})]\vb_n|\nonumber\\
[-8pt]\\ [-8pt]
&&\qquad=O_p\bigl(\sqrt{np_n}\bigr),\nonumber
\end{eqnarray}
\begin{eqnarray}\label{cat2}
&&\sup_{\Vert\vbeta_n-\vbeta_{n0}\Vert\leq \Delta\sqrt{p_n/n}}
\sup_{\Vert\vb_n\Vert=1}|\vb_n^T[\vE_n(\vbeta_n)-\overline{\vE}_n(\vbeta_{n})]\vb_n|\nonumber\\
[-8pt]\\ [-8pt]
&&\qquad=O_p\bigl(\sqrt{np_n}\bigr),\nonumber
\end{eqnarray}
\begin{eqnarray}\label{cat3}
&&\sup_{\Vert\vbeta_n-\vbeta_{n0}\Vert\leq \Delta\sqrt{p_n/n}}
\sup_{\Vert\vb_n\Vert=1}|\vb_n^T[\vG_n(\vbeta_n)-\overline{\vG}_n(\vbeta_{n})]\vb_n|\nonumber\\
[-8pt]\\ [-8pt]
&&\qquad=O_p\bigl(\sqrt{np_n}\bigr).\nonumber
\end{eqnarray}
We have
\begin{eqnarray*}
&&|\vb_n^T[\vH_n(\vbeta_n)-\overline{\vH}_n(\vbeta_{n})]\vb_n|\\
&&\qquad=\Biggl|\sum_{i=1}^n\vb_n^T\vX_i^T\vA_i^{1/2}(\vbeta_n)
[\widehat{\vR}^{-1}-\overline{\vR}^{-1}]\vA_i^{1/2}(\vbeta_{n})\vX_i\vb_n\Biggr|\\
&&\qquad\leq\Vert\widehat{\vR}^{-1}-\overline{\vR}^{-1}\Vert\lambda_{\max}(\vA_i(\vbeta_{n}))
\lambda_{\max}\Biggl(\sum_{i=1}^n\vX_i^T\vX_i\Biggr)\Vert\vb_n\Vert^2.
\end{eqnarray*}
By assumptions (A2) and
(A4), (\ref{cat1}) is proved. Next, note that
\begin{eqnarray*}
&&|\vb_n^T[\vE_n(\vbeta_n)-\overline{\vE}_n(\vbeta_{n})]\vb_n|\\
&&\qquad=
\frac{1}{2}\Biggl|\sum_{i=1}^n\sum_{j=1}^m\bigl(1-2\pi_{ij}(\vbeta_n)\bigr)\epsilon_{ij}(\vbeta_{n})
\vb_n^T\vX_i^T\vA_i^{1/2}(\vbeta_n)\\
&&\hspace*{127pt}{}\times[\widehat{\vR}^{-1}-\overline{\vR}^{-1}]\ve_j\ve_j^T\vX_i\vb_n\Biggr|\\
&&\qquad\leq
\sum_{i=1}^n\sum_{j=1}^m\vA_{ij}^{-1/2}(\vbeta_n)|\vb_n^T\vX_i^T\vA_i^{1/2}(\vbeta_n)
[\widehat{\vR}^{-1}-\overline{\vR}^{-1}]\ve_j| \cdot
|\ve_j^T\vX_i\vb_n|\\
& &\qquad\leq \sum_{i=1}^n\sum_{j=1}^m\vA_{ij}^{-1/2}(\vbeta_n)
\Vert\widehat{\vR}^{-1}-\overline{\vR}^{-1}\Vert\cdot
\Vert\vA_{i}^{1/2}(\vbeta_n)\Vert\cdot \Vert\vX_i\vb_n\Vert^2.
\end{eqnarray*}
Thus,
\begin{eqnarray*}
&&\sup_{\Vert\vbeta_n-\vbeta_{n0}\Vert\leq
\Delta\sqrt{p_n/n}}\sup_{\Vert\vb_n\Vert=1}|\vb_n^T[\vE_n(\vbeta_n)-\overline{\vE}_n(\vbeta_{n})]\vb_n|\\
&&\qquad\leq C\Vert\widehat{\vR}^{-1}-\overline{\vR}^{-1}\Vert\cdot
\sum_{i=1}^n\sum_{j=1}^m \sup_{\Vert\vbeta_n-\vbeta_{n0}\Vert\leq
\Delta\sqrt{p_n/n}}\vA_{ij}^{-1/2}(\vbeta_n)
\sup_{\Vert\vb_n\Vert=1}\Vert\vX_i\vb_n\Vert^2\\
&&\qquad=O_p\bigl(\sqrt{p_n/n}\bigr)O(n)=O_p\bigl(\sqrt{np_n}\bigr),
\end{eqnarray*}
by assumption (A3).
(\ref{cat3}) is proved
similarly.
\end{pf*}

\begin{pf*}{Proof of Lemma \ref{tiger4}} By (\ref{rep}), it is sufficient to
verify that
\begin{eqnarray}\label{Enbn}
\sup_{\Vert\vbeta_n-\vbeta_{n0}\Vert\leq \Delta\sqrt{p_n/n}}\
\sup_{\Vert\vb_n\Vert=1}|\vb_n^T\overline{\vE}_n(\vbeta_n)\vb_n|
&=&O_p\bigl(\sqrt{n}p_n\bigr),\\ \label{Gnbn}
\sup_{\Vert\vbeta_n-\vbeta_{n0}\Vert\leq
\Delta\sqrt{p_n/n}}\
\sup_{\Vert\vb_n\Vert=1}|\vb_n^T\overline{\vG}_n(\vbeta_n)\vb_n|
&=&O_p\bigl(\sqrt{n}p_n\bigr).
\end{eqnarray}
First, note that we have the following
decomposition of $\overline{\vE}_n(\vbeta_n)$:
\begin{eqnarray*}
\overline{\vE}_n(\vbeta_n)&=&\frac{1}{2}\sum_{i=1}^n\sum_{j=1}^m\bigl(1-2\pi_{ij}(\vbeta_{n0})\bigr)\epsilon_{ij}(\vbeta_{n0})
\vX_i^T\vA_i^{1/2}(\vbeta_{n0})\overline{\vR}^{-1}\ve_j\ve_j^T\vX_i\\
&&{}+\frac{1}{2}\sum_{i=1}^n\sum_{j=1}^m\bigl(1-2\pi_{ij}(\vbeta_{n0})\bigr)\epsilon_{ij}(\vbeta_{n0})\vX_i^T[\vA_i^{1/2}(\vbeta_{n})-\vA_i^{1/2}(\vbeta_{n0})]\\
&&\hphantom{{}+\frac{1}{2}\sum_{i=1}^n\sum_{j=1}^m}{}\times\overline{\vR}^{-1}\ve_j\ve_j^T\vX_i\\
&&{}+\frac{1}{2}\sum_{i=1}^n\sum_{j=1}^m\bigl[\bigl(1-2\pi_{ij}(\vbeta_{n})\bigr)\vA_{ij}^{-1/2}(\vbeta_{n})-
\bigl(1-2\pi_{ij}(\vbeta_{n0})\bigr)\vA_{ij}^{-1/2}(\vbeta_{n0})\bigr]
\\&&\hphantom{{}+\frac{1}{2}\sum_{i=1}^n\sum_{j=1}^m}{}\times\bigl(Y_{ij}-\pi_{ij}(\vbeta_{n0})\bigr)\vX_i^T\vA_i^{1/2}(\vbeta_{n})\overline{\vR}^{-1}\ve_j\ve_j^T\vX_i\\
&&{}+\frac{1}{2}\sum_{i=1}^n\sum_{j=1}^m\bigl(1-2\pi_{ij}(\vbeta_{n})\bigr)\vA_{ij}^{-1/2}(\vbeta_{n})
\bigl(\pi_{ij}(\vbeta_{n0})-\pi_{ij}(\vbeta_{n})\bigr)\\
&&\hphantom{{}+\frac{1}{2}\sum_{i=1}^n\sum_{j=1}^m}{}\times\vX_i^T\vA_i^{1/2}(\vbeta_{n})\overline{\vR}^{-1}\ve_j\ve_j^T\vX_i\\
&\triangleq & \overline{\vE}_{1n}(\vbeta_{n0})+
\sum_{k=2}^4\overline{\vE}_{kn}(\vbeta_{n}).
\end{eqnarray*}
Thus, to prove (\ref{Enbn}),
it suffices to verify that $\sup_{\Vert\vb_n\Vert=1}|\vb_n^T\overline{\vE}_{1n}(\vbeta_{n0})\vb_n|=O_P(\sqrt{n}p_n)$
and $\sup_{\Vert\vbeta_n-\vbeta_{n0}\Vert\leq \Delta\sqrt{p_n/n}}\
\sup_{\Vert\vb_n\Vert=1}|\vb_n^T\overline{\vE}_{kn}(\vbeta_n)\vb_n|=O_P(\sqrt{n}p_n)$.
We first prove that
$\sup_{\Vert\vb_n\Vert=1}|\vb_n^T\overline{\vE}_{1n}(\vbeta_{n0})\vb_n|=O_P(\sqrt{n}p_n),$
by verifying that $\Vert\overline{\vE}_{1n}(\vbeta_{n0})\Vert=O_P(\sqrt{n}p_n)$,
where
$\Vert\overline{\vE}_{1n}(\vbeta_{n0})\Vert=\sqrt{\operatorname{trace}(\overline{\vE}_{1n}(\vbeta_{n0})
\overline{\vE}_{1n}^T(\vbeta_{n0}))}$:
\begin{eqnarray*}
&&E[\Vert\overline{\vE}_{1n}(\vbeta_{n0})\Vert^2]\\
&&\qquad= \frac{1}{4}\sum_{i=1}^n\sum_{j_1=1}^m\sum_{j_2=1}^m
\bigl(1-2\pi_{ij_1}(\vbeta_{n0})\bigr)\bigl(1-2\pi_{ij_2}(\vbeta_{n0})\bigr)
E[\epsilon_{ij_1}(\vbeta_{n0})\epsilon_{ij_2}(\vbeta_{n0})]\\
&&\qquad\hphantom{= \frac{1}{4}\sum_{i=1}^n\sum_{j_1=1}^m\sum_{j_2=1}^m}
{}\times\operatorname{trace}[\vX_i^T\vA_i^{1/2}(\vbeta_{n0})\overline{\vR}^{-1}\ve_{j_1}\ve_{j_1}^T\vX_i
\vX_i^T\ve_{j_2}\ve_{j_2}^T\\
&&\qquad\hspace*{179pt}{}\times\overline{\vR}^{-1}\vA_i^{1/2}(\vbeta_{n0})\vX_i]\\
&&\qquad\leq  C\sum_{i=1}^n\sum_{j_1=1}^m\sum_{j_2=1}^m |\ve_{j_1}^T\vX_i
\vX_i^T\ve_{j_2}\ve_{j_2}^T\overline{\vR}^{-1}\vA_i^{1/2}(\vbeta_{n0})\\
&&\qquad\hspace*{90pt}{}\times\vX_i
\vX_i^T\vA_i^{1/2}(\vbeta_{n0})\overline{\vR}^{-1}\ve_{j_1}|\\
&&\qquad\leq  C\sum_{i=1}^n\sum_{j_1=1}^m\sum_{j_2=1}^m \Vert\ve_{j_1}^T\vX_i\Vert\cdot
\Vert\vX_i^T\ve_{j_2}\Vert\cdot
\Vert\ve_{j_2}^T\overline{\vR}^{-1}\vA_i^{1/2}(\vbeta_{n0})\vX_i\Vert\\
&&\qquad\hphantom{\leq  C\sum_{i=1}^n\sum_{j_1=1}^m\sum_{j_2=1}^m}{}\times \Vert\vX_i^T\vA_i^{1/2}(\vbeta_{n0})\overline{\vR}^{-1}\ve_{j_1}\Vert.
\end{eqnarray*}
Note that $\Vert\ve_{j_1}^T\vX_i\Vert=\Vert\vX_{ij_1}\Vert$,
$\Vert\vX_i^T\ve_{j_2}\Vert=\Vert\vX_{ij_2}\Vert$,
$\Vert\ve_{j_2}^T\overline{\vR}^{-1}\vA_i^{1/2}(\vbeta_{n0})\vX_i\Vert\leq\break
C(\operatorname{trace}(\vX_i\vX_i^T))^{1/2}$ and
$\Vert\vX_i^T\vA_i^{1/2}(\vbeta_{n0})\overline{\vR}^{-1}\ve_{j_1}\Vert\leq
C(\operatorname{trace}(\vX_i\vX_i^T))^{1/2}$. Thus,\vadjust{\eject}
\begin{eqnarray*}
E[\Vert\overline{\vE}_{1n}(\vbeta_{n0})\Vert^2]
&\leq &   C\sum_{i=1}^n\sum_{j_1=1}^m\sum_{j_2=1}^m\Vert\vX_{ij_1}\Vert \cdot
\Vert\vX_{ij_2}\Vert\operatorname{trace}(\vX_i\vX_i^T)\\
&\leq & C \cdot \max_{i,j}\Vert\vX_{ij}\Vert^2
\operatorname{trace}\Biggl(\sum_{i=1}^n\vX_i\vX_i^T\Biggr)=O(np_n^2),
\end{eqnarray*}
by assumptions (A1) and (A3). This
implies that
$\sup_{\Vert\vb_n\Vert=1}|\vb_n^T\overline{\vE}_{1n}(\vbeta_{n0})\vb_n|=O_p(\sqrt{n}p_n)$.
Next, we have
\begin{eqnarray*}
&& |\vb_n^T\overline{\vE}_{2n}(\vbeta_n)\vb_n|\\
&&\qquad=\Biggl|\frac{1}{2}\sum_{i=1}^n\sum_{j=1}^m\bigl(1-2\pi_{ij}(\vbeta_{n0})\bigr)\epsilon_{ij}^{1/2}(\vbeta_{n0})
\vb_n^T\vX_i^T[\vA_i^{1/2}(\vbeta_{n})-\vA_i^{1/2}(\vbeta_{n0})]\\
&&\qquad\hspace*{225pt}{}\times\overline{\vR}^{-1}\ve_j\ve_j^T\vX_i\vb_n\Biggr|\\
&&\qquad\leq  C\sum_{i=1}^n\sum_{j=1}^m
|\vb_n^T\vX_i^T[\vA_i^{1/2}(\vbeta_{n})-\vA_i^{1/2}(\vbeta_{n0})]\overline{\vR}^{-1}\ve_j|\cdot
|\ve_j^T\vX_i\vb_n|\\
& &\qquad\leq  C\sum_{i=1}^n\sum_{j=1}^m
\Vert\vX_i\vb_n\Vert^2\lambda_{\max}(\overline{\vR}^{-1})\max_{j}|\vA_{ij}^{1/2}(\vbeta_{n})-\vA_{ij}^{1/2}(\vbeta_{n0})|.
\end{eqnarray*}
Note that there exists some $\vbeta_{n}^*$ between $\vbeta_{n}$ and
$\vbeta_{n0}$ such that
\begin{eqnarray*}
\vA_{ij}^{1/2}(\vbeta_{n})-\vA_{ij}^{1/2}(\vbeta_{n0})&=&
\tfrac{1}{2}\vA_{ij}^{1/2}(\vbeta_{n}^*)\bigl(1-2\pi_{ij}(\vbeta_{n}^*)\bigr)\vX_{ij}^T(\vbeta_{n}-\vbeta_{n0})\\
&\leq & C \Vert\vX_{ij}\Vert \cdot \Vert\vbeta_{n}-\vbeta_{n0}\Vert.
\end{eqnarray*}
Therefore,
\begin{eqnarray*}
&&\sup_{\Vert\vbeta_n-\vbeta_{n0}\Vert\leq
\Delta\sqrt{p_n/n}}\sup_{\Vert\vb_n\Vert=1}|\vb_n^T\overline{\vE}_{2n}(\vbeta_n)\vb_n|\\
&&\qquad\leq C \max_{i,j}\Vert\vX_{ij}\Vert \sup_{\Vert\vbeta_n-\vbeta_{n0}\Vert\leq
\Delta\sqrt{p_n/n}}\Vert\vbeta_n-\vbeta_{n0}\Vert\cdot
\lambda_{\max}\Biggl(\sum_{i=1}^n\vX_i^T\vX_i\Biggr) \\
&&\qquad=O\bigl(\sqrt{p_n}\bigr)O\bigl(\sqrt{p_n/n}\bigr)O(n)=O\bigl(\sqrt{n}p_n\bigr).
\end{eqnarray*}
Similarly, we can show
that $\sup_{\Vert\vbeta_n-\vbeta_{n0}\Vert\leq
\Delta\sqrt{p_n/n}}\sup_{\Vert\vb_n\Vert=1}|\vb_n^T\overline{\vE}_{kn}(\vbeta_n)\vb_n|=\break O(\sqrt{n}p_n)$,
$k=3,4$. This proves (\ref{Enbn}). Similarly, we can
prove (\ref{Gnbn}).
\end{pf*}

\begin{pf*}{Proof of Lemma \ref{tiger5}}
The proof is given in the online supplementary material.
\end{pf*}

\begin{pf*}{Proof of Lemma \ref{tiger6}}
We write
$\valpha_n^T\overline{\vM}^{-1/2}_n(\vbeta_{n0})\overline{\vS}_n(\vbeta_{n0})=\sum_{i=1}^nZ_{ni}$,
where $
Z_{ni}=\valpha_n^T\overline{\vM}^{-1/2}_n(\vbeta_{n0})\vX_i^T\vA_i^{1/2}(\vbeta_{n0})
\overline{\vR}^{-1}\vepsilon_i(\vbeta_{n0}).$  Since
$\overline{\vM}_n(\vbeta_{n0})=\break \Cov(\overline{\vS}_n(\vbeta_{n0}))$, we have
$\Var(\valpha_n^T\overline{\vM}^{-1/2}_n(\vbeta_{n0})\overline{\vS}_n(\vbeta_{n0}))=1$.
To establish the asymptotic normality, it suffices to check the
Lindberg condition, that is, $\forall \epsilon>0$, $
\sum_{i=1}^nE[Z_{ni}^2I(|Z_{ni}|>\epsilon)]\ra 0.$ By
the Cauchy--Schwarz inequality,
\begin{eqnarray*}
Z^2_{ni}&\leq&
\Vert\valpha_n^T\overline{\vM}^{-1/2}_n(\vbeta_{n0})\vX_i^T\vA_i^{1/2}(\vbeta_{n0})
\overline{\vR}^{-1}\Vert^2\cdot \Vert\vepsilon_i(\vbeta_{n0})\Vert^2\\
&\leq& \lambda_{\max}(\overline{\vR}^{-2})\lambda_{\max}(\vA_i(\vbeta_{n0}))
(\valpha_n^T\overline{\vM}^{-1/2}_n(\vbeta_{n0})\vX_i^T\vX_i\overline{\vM}^{-1/2}_n(\vbeta_{n0})\valpha_n)\\
&&\times\Vert\vepsilon_i(\vbeta_{n0})\Vert^2\\
&\leq & C\gamma_{ni}\Vert\vepsilon_i(\vbeta_{n0})\Vert^2,
\end{eqnarray*}
where
$\gamma_{ni}\triangleq
\valpha_n^T\overline{\vM}^{-1/2}_n(\vbeta_{n0})\vX_i^T\vX_i\overline{\vM}^{-1/2}_n(\vbeta_{n0})\valpha_n$.
Next, we will show that\break $\max_{1\leq i \leq n}\gamma_{ni}\ra 0$ as $n\ra
\infty$. Note that $ \gamma_{ni}\leq
\lambda_{\max}(\vX_i^T\vX_i)\lambda_{\min}^{-1}(\overline{\vM}_n(\vbeta_{n0})).
$ Since $\overline{\vM}_n(\vbeta_{n0})$
is symmetric, to evaluate $\lambda_{\min}(\overline{\vM}_n(\vbeta_{n0}))$, $\forall \vb_n\in R^{p_n}$, we have
\begin{eqnarray*}
\vb_n^T\overline{\vM}_n(\vbeta_{n0})\vb_n&\geq&
\lambda_{\min}(\vR_0)\lambda_{\min}(\overline{\vR}^{-2})\sum_{i=1}^n\lambda_{\min}(\vA_i(\vbeta_{n0}))
\vb_n^T\vX_i^T\vX_i\vb_n\\
&\geq & C\vb_n^T\Biggl(\sum_{i=1}^n\vX_i^T\vX_i\Biggr)\vb_n\geq C
\lambda_{\min}\Biggl(\sum_{i=1}^n\vX_i^T\vX_i\Biggr)\Vert b_n\Vert^2.
\end{eqnarray*}
Thus, $\inf_{\Vert\vb_n\Vert=1}|\vb_n^T\overline{\vM}_n(\vbeta_{n0})\vb_n|\geq C
\lambda_{\min}(\sum_{i=1}^n\vX_i^T\vX_i)$ and this
implies that $\lambda_{\min}
(\overline{\vM}_n(\vbeta_{n0}))\geq C
\lambda_{\min}(\sum_{i=1}^n\vX_i^T\vX_i)$. Therefore, we have
\[
\gamma_{ni}\leq
\frac{\lambda_{\max}(\vX_i^T\vX_i)}{C\lambda_{\min}(\sum_{i=1}^n\vX_i^T\vX_i)}
\leq \frac{\operatorname{Tr}(\vX_i^T\vX_i)}{C\lambda_{\min}(\sum_{i=1}^n\vX_i^T\vX_i)}
=\frac{\sum_{j=1}^m\vX_{ij}^T\vX_{ij}}{C\lambda_{\min}(\sum_{i=1}^n\vX_i^T\vX_i)}.
\]
It follows that $\max_{1\leq i \leq n}\gamma_{ni}\leq O(n^{-1}p_n)=o(1)$. We have
\[
\sum_{i=1}^nE[Z_{ni}^2I(|Z_{ni}|>\epsilon)] \leq \sum_{i=1}^n
C\gamma_{ni} E\biggl[\Vert\vepsilon_i(\vbeta_{n0})\Vert^2
I\biggl\{\Vert\vepsilon_i(\vbeta_{n0})\Vert^2>\frac{\epsilon^2}{C\gamma_{ni}}\biggr\}\biggr].
\]
Note that $\Vert\vepsilon_i(\vbeta_{n0})\Vert^2$ is uniformly bounded, by
assumption (A2). Thus, for all $\epsilon>0$ and $\delta>0$, there exists  a positive integer $N$ such that (1)
$I\{\Vert\vepsilon_i(\vbeta_{n0})\Vert^2>\frac{\epsilon^2}{C\gamma_{ni}}\}=0$
for all  $n>N$; (2) $\sum_{i=1}^N C\gamma_{ni}\leq \delta$ for all $n$
sufficiently large. This ensures that
\[
\sum_{i=1}^n C\gamma_{ni}
E\biggl[\Vert\vepsilon_i(\vbeta_{n0})\Vert^2
I\biggl\{\Vert\vepsilon_i(\vbeta_{n0})\Vert^2>\frac{\epsilon^2}{C\gamma_{ni}}\biggr\}\biggr]\ra 0.
\]
Therefore, the Lindberg condition is verified.
\end{pf*}

\begin{pf*}{Proof of Theorem \ref{cat5}} It is sufficient to show that for
$\vb_n\in R^{p_n}$,
\begin{equation}\label{lion2}
\sup_{\Vert\vb_n\Vert=1}|\vb_n^T(\widehat{\vSigma}_n-\vSigma_n)\vb_n|=o_p(n^{-1}).
\end{equation}
We use the
conclusion of Theorem \ref{thm1} throughout the proof. Note that we can write
$\widehat{\vSigma}_n-\vSigma_n=I_{n1}+I_{n2}+I_{n3}$, where
\begin{eqnarray*}
I_{n1}&=&\vH^{-1}_n(\widehat{\vbeta}_{n})
[\widehat{\vM}_n(\widehat{\vbeta}_{n})-\overline{\vM}_n(\vbeta_{n0})]
\vH^{-1}_n(\widehat{\vbeta}_{n}),\\
I_{n2}&=&[\vH^{-1}_n(\widehat{\vbeta}_{n})-\overline{\vH}^{-1}_n(\vbeta_{n0})]
\overline{\vM}_n(\vbeta_{n0})
\vH^{-1}_n(\widehat{\vbeta}_{n}),\\
I_{n3}&=& \overline{\vH}^{-1}_n(\vbeta_{n0})\overline{\vM}_n(\vbeta_{n0})
[\vH^{-1}_n(\widehat{\vbeta}_{n})-\overline{\vH}^{-1}_n(\vbeta_{n0})].
\end{eqnarray*}
Thus, (\ref{lion2}) is implied by $\sup_{\Vert\vb_n\Vert=1}|\vb_n^TI_{ni}\vb_n|=o_p(1)$. We have
\begin{eqnarray*}
&&\sup_{\Vert\vb_n\Vert=1}|\vb_n^TI_{n1}\vb_n|\\
&&\qquad\leq \frac{\max(|\lambda_{\max}(\widehat{\vM}_n(\widehat{\vbeta}_{n})-\overline{\vM}_n(\vbeta_{n0}))|,
|\lambda_{\min}(\widehat{\vM}_n(\widehat{\vbeta}_{n})-\overline{\vM}_n(\vbeta_{n0}))|)}
{\lambda_{\min}^2(\vH_n(\widehat{\vbeta}_{n}))}.
\end{eqnarray*}
To evaluate the
eigenvalues of
$\widehat{\vM}_n(\widehat{\vbeta}_{n})-\overline{\vM}_n(\vbeta_{n0})$, we have
\begin{eqnarray*}
&&|\vc_n^T[\widehat{\vM}_n(\widehat{\vbeta}_{n})-\overline{\vM}_n(\vbeta_{n0})]\vc_n|\\
&&\qquad\leq
|\vc_n^T[\widehat{\vM}_n(\widehat{\vbeta}_{n})-\widehat{\vM}_n(\vbeta_{n0})]\vc_n|+
|\vc_n^T[\widehat{\vM}_n(\vbeta_{n0})-\overline{\vM}_n(\vbeta_{n0})]\vc_n|
\end{eqnarray*}
for
$\vc_n\in R^{p_n}$. Note that
\begin{eqnarray*}
&&\sup_{\Vert\vc_n\Vert=1}|\vc_n^T[\widehat{\vM}_n(\widehat{\vbeta}_{n})-\widehat{\vM}_n(\vbeta_{n0})]\vc_n|\\
&&\qquad\leq
\sup_{\Vert\vc_n\Vert=1}\Biggl|\sum_{i=1}^n\vc_n^T\vX_i^T[\vA_i^{1/2}(\widehat{\vbeta}_{n})-\vA_i^{1/2}(\vbeta_{n0})]
\widehat{\vR}^{-1}\vepsilon_i(\widehat{\vbeta}_{n})\vepsilon_i^T(\widehat{\vbeta}_{n})\\
&&\qquad\hspace*{175pt}{}\times\widehat{\vR}^{-1}\vA_i^{1/2}(\widehat{\vbeta}_{n})\vX_i\vc_n\Biggr|\\
&&\qquad\quad{}+ \sup_{\Vert\vc_n\Vert=1}\Biggl|\sum_{i=1}^n\vc_n^T\vX_i^T\vA_i^{1/2}(\vbeta_{n0})
\widehat{\vR}^{-1}\vepsilon_i(\widehat{\vbeta}_{n})\vepsilon_i^T(\widehat{\vbeta}_{n})
\widehat{\vR}^{-1}\\
&&\qquad\hspace*{100pt}{}\times[\vA_i^{1/2}(\widehat{\vbeta}_{n})-\vA_i^{1/2}(\vbeta_{n0})]\vX_i\vc_n\Biggr|\\
&&\qquad\quad{}+\sup_{\Vert\vc_n\Vert=1}\Biggl|\sum_{i=1}^n\vc_n^T\vX_i^T\vA_i^{1/2}(\vbeta_{n0})
\widehat{\vR}^{-1}[\vepsilon_i(\widehat{\vbeta}_{n})\vepsilon_i^T(\widehat{\vbeta}_{n})-
\vepsilon_i(\vbeta_{n0})\vepsilon_i^T(\vbeta_{n0})]\\
&&\qquad\hspace*{215pt}{}\times\widehat{\vR}^{-1}\vA_i^{1/2}(\vbeta_{n0})\vX_i\vc_n\Biggr|\\
&&\qquad\triangleq  \sup_{\Vert\vc_n\Vert=1}J_{n1}+\sup_{\Vert\vc_n\Vert=1}J_{n2}+\sup_{\Vert\vc_n\Vert=1}J_{n3}.
\end{eqnarray*}
Note that
\[
J_{n1} \leq
\sum_{i=1}^n\Vert\vc_n^T\vX_i^T[\vA_i^{1/2}(\widehat{\vbeta}_{n})-\vA_i^{1/2}(\vbeta_{n0})]\Vert\cdot
\Vert\widehat{\vR}^{-1}\vepsilon_i(\widehat{\vbeta}_{n})\Vert^2\cdot
\Vert\vA_i^{1/2}(\widehat{\vbeta}_{n})\vX_i\vc_n\Vert.
\]
We have
$\Vert\vA_i^{1/2}(\widehat{\vbeta}_{n})\vX_i\vc_n\Vert\leq \Vert\vX_i\vc_n\Vert$ and
\begin{eqnarray*}
\Vert\vc_n^T\vX_i^T[\vA_i^{1/2}(\widehat{\vbeta}_{n})-\vA_i^{1/2}(\vbeta_{n0})]\Vert
&\leq&
\Vert\vX_i\vc_n\Vert\max_{j}|\vA_{ij}^{1/2}(\widehat{\vbeta}_{n})-\vA_{ij}^{1/2}(\vbeta_{n0})|\\
&\leq & C \Vert\vX_i\vc_n\Vert \cdot \Vert\vX_{ij}\Vert \cdot
\Vert\widehat{\vbeta}_{n}-\vbeta_{n0}\Vert.
\end{eqnarray*}
Furthermore,
\begin{eqnarray*}
\Vert\widehat{\vR}^{-1}\vepsilon_i(\widehat{\vbeta}_{n})\Vert^2&=&\bigl(\vY_i-\vpi_i(\widehat{\vbeta}_{n})\bigr)^T\vA_i^{-1/2}(\widehat{\vbeta}_{n})\widehat{\vR}^{-2}
\vA_i^{-1/2}(\widehat{\vbeta}_{n})\bigl(\vY_i-\vpi_i(\widehat{\vbeta}_{n})\bigr)\\
&\leq&\lambda_{\max}(\widehat{\vR}^{-2})\lambda_{\max}(\vA_i^{-1}(\widehat{\vbeta}_{n}))
\Vert\vY_i-\vpi_i(\widehat{\vbeta}_{n})\Vert^2\leq CO_p(1).
\end{eqnarray*}
Thus,
\[
\sup_{\Vert\vc_n\Vert=1}J_{n1}\leq
O_p(1)\Vert\widehat{\vbeta}_{n}-\vbeta_{n0}\Vert\max_{i,j}\Vert\vX_{ij}\Vert
\lambda_{\max}\Biggl(\sum_{i=1}^n\vX_i^T\vX_i\Biggr)=o_p(n).
\]
Similarly, $\sup_{\Vert\vc_n\Vert=1}J_{n2}=o_p(n)$ and $\sup_{\Vert\vc_n\Vert=1}J_{n3}=o_p(n)$. Thus,
\[
\sup_{\Vert\vc_n\Vert=1}|\vc_n^T[\widehat{\vM}_n(\widehat{\vbeta}_{n})-\widehat{\vM}_n(\vbeta_{n0})]\vc_n|=o_p(n).
\]
Similarly,
$\sup_{\Vert\vc_n\Vert=1}|\vc_n^T[\widehat{\vM}_n(\vbeta_{n0})-\overline{\vM}_n(\vbeta_{n0})]\vc_n|=o_p(n)$.
Finally, note that
\begin{eqnarray*}
\lambda_{\min}(\vH_n(\widehat{\vbeta}_{n}))&\geq&
\lambda_{\min}(\widehat{\vR})\min_{i,j}\bigl(\pi_{ij}(\widehat{\vbeta}_{n})
\bigl(1-\pi_{ij}(\widehat{\vbeta}_{n})\bigr)\bigr)\lambda_{\min}
\Biggl(\sum_{i=1}^n\vX_i^T\vX_i\Biggr)\\
&=&O_p(n).
\end{eqnarray*}
Thus,
$\sup_{\Vert\vb_n\Vert=1}|\vb_n^TI_{n1}\vb_n|=o_p(n^{-1})$. We can also prove that
$\sup_{\Vert\vb_n\Vert=1}|\vb_n^TI_{ni}\times\break\vb_n|=o_p(n^{-1})$, $i=2,3$, by first noting that
\[
\vH^{-1}_n(\widehat{\vbeta}_{n})-\overline{\vH}^{-1}_n(\vbeta_{n0})=
[\vH^{-1}_n(\widehat{\vbeta}_{n})-\overline{\vH}^{-1}_n(\widehat{\vbeta}_{n})]+
[\overline{\vH}^{-1}_n(\widehat{\vbeta}_{n})-\overline{\vH}^{-1}_n(\vbeta_{n0})]
\]
and then using the expressions
\begin{eqnarray*}
\vH^{-1}_n(\widehat{\vbeta}_{n})-\overline{\vH}^{-1}_n(\widehat{\vbeta}_{n})
&=&\overline{\vH}^{-1}_n(\widehat{\vbeta}_{n})[\overline{\vH}_n(\widehat{\vbeta}_{n})-
\vH_n(\widehat{\vbeta}_{n})]\vH^{-1}_n(\widehat{\vbeta}_{n}),\\
\overline{\vH}^{-1}_n(\widehat{\vbeta}_{n})-\overline{\vH}^{-1}_n(\vbeta_{n0})
&=&\overline{\vH}^{-1}_n(\vbeta_{n0})[\overline{\vH}_n(\vbeta_{n0})-
\overline{\vH}_n(\widehat{\vbeta}_{n})]\overline{\vH}^{-1}_n(\widehat{\vbeta}_{n}).
\end{eqnarray*}
\upqed\end{pf*}

\begin{pf*}{Proof of Corollary \ref{wald}} It is sufficient to show that
%
\begin{equation}\label{dog3}
[(\vL_n\widehat{\vSigma}_n\vL_n^T)^{-1/2}-(\vL_n\vSigma_n\vL_n^T)^{-1/2}]\vL_n(\widehat{\vbeta}_n-\vbeta_{n0})\ra 0
\end{equation}
in probability. Note that the left-hand side can be written as
\[
[(\vL_n\widehat{\vSigma}_n\vL_n^T)^{-1/2}(\vL_n\vSigma_n\vL_n^T)^{1/2}-\vI_l]
(\vL_n\vSigma_n\vL_n^T)^{-1/2}\vL_n(\widehat{\vbeta}_n-\vbeta_{n0})
\]
and thus
(\ref{dog3}) is implied by
\[
(\vL_n\widehat{\vSigma}_n\vL_n^T)^{-1}(\vL_n\vSigma_n\vL_n^T)-\vI_l
=(\vL_n\widehat{\vSigma}_n\vL_n^T)^{-1}\vL_n(\vSigma_n-\widehat{\vSigma}_n)\vL_n^T=o_p(1).
\]
Let $\vu_i$ denote the $l\times 1$ unit vector with the $i$th element
being 1 and all of the other elements being 0. Then, for all $1\leq i,j\leq
l$, we have, by the Cauchy--Schwarz inequality,
\begin{eqnarray*}
&&|\vu_i^T(\vL_n\widehat{\vSigma}_n\vL_n^T)^{-1}\vL_n(\vSigma_n-\widehat{\vSigma}_n)\vL_n^T\vu_j|\\
& &\qquad\leq |\vu_i^T(\vL_n\widehat{\vSigma}_n\vL_n^T)^{-2}\vu_i|^{1/2}
|\vu_j^T[\vL_n(\vSigma_n-\widehat{\vSigma}_n)\vL_n^T]^2\vu_j|^{1/2}\\
&&\qquad\leq
\frac{\Vert\vL_n(\vSigma_n-\widehat{\vSigma}_n)\vL_n^T\Vert}{\lambda_{\min}(\vL_n\widehat{\vSigma}_n\vL_n^T)}.
\end{eqnarray*}
Now, for any $l$-dimensional vector such that $\Vert\vb\Vert=1$, we have
\begin{eqnarray*}
|\vb^T\vL_n\widehat{\vSigma}_n\vL_n^T\vb| &\geq&
|\vb^T\vL_n\vSigma_n\vL_n^T\vb|-
|\vb^T\vL_n(\widehat{\vSigma}_n-\vSigma_n)\vL_n^T\vb|\\
&\geq& \lambda_{\min}(\vSigma_n)+o_p(n^{-1})\\
&\geq& \frac{\lambda_{\min}(\overline{\vM}_n(\vbeta_{n0}))}
{\lambda^2_{\max}(\overline{\vH}_n(\vbeta_{n0}))}+o_p(n^{-1}),
\end{eqnarray*}
where the
second inequality uses Theorem \ref{cat5}. By (\ref{mn}),
$\lambda_{\min}(\overline{\vM}_n(\vbeta_{n0}))\geq
c_1\lambda_{\min}(\sum_{i=1}^n\vX_i^T\vX_i)$ for some positive constant $c_1$.
Similarly, we can show that $\lambda_{\max}(\overline{\vH}_n(\vbeta_{n0}))\leq
c_2\lambda_{\max}(\sum_{i=1}^n\vX_i^T\vX_i)$ for some positive constant $c_2$.
Thus, $\lambda_{\min}(\vL_n\widehat{\vSigma}_n\vL_n^T)\geq O_p(n^{-1})$. This
proves that
$\frac{\Vert\vL_n(\vSigma_n-\widehat{\vSigma}_n)\vL_n^T\Vert}{\lambda_{\min}(\vL_n\widehat{\vSigma}_n\vL_n^T)}=o_p(1),$
by Theorem \ref{cat5}.
\end{pf*}
\end{appendix}

\section*{Acknowledgments}
The author would like to thank the Associate Editor and two referees for their
constructive and insightful comments that significantly improved this paper.

\begin{supplement}[id=suppA]
\stitle{Supplement to ``GEE analysis of clustered binary data with
diverging number of covariates''}
\sfilename{supplement.pdf}
\slink[doi]{10.1214/10-AOS846SUPP}  
\sdatatype{.pdf}
\sdescription{The proofs of (\ref{inia}), Lemma \ref{tiger5}, (\ref{high}) and
Theorem 5.1 are provided in this supplementary article [Wang~(\citeyear{wan2010})].}
\end{supplement}

\printaddresses

\end{document}